\magnification=\magstep1 \hoffset .65pt
\def\R{I\!\! R}

\def\C{C \!\! \! \!  I\, }

\def\RP{I\!\! R I\!\! P}
\def\Ri{Riemannian }

\input epsf.tex

\parindent=8pt

\parindent=8pt
\hskip11cm {\bf DRAFT. April 24, 2004}
\bigskip
\centerline{\bf  Fitting Hyperbolic pants
to a three-body problem. }
\bigskip
\bigskip

{\bf  Abstract.}  Consider the three-body problem
 with an  attractive $1/r^2$
potential.  Modulo symmetries,
the dynamics of the    bounded zero-angular momentum solutions  
is equivalent to a  geodesic flow on the  thrice-punctured sphere,
or ``pair of pants''.  The sphere is the shape sphere.
  The punctures are the binary collisions.  The metric generating
the geodesics is the Jacobi-Maupertuis metric.  The metric is 
complete, has infinite area, and its ends, the neighborhoods
of the punctures,  are asymptotically cylindrical.
Our main result is  that when the three masses are  equal then the
metric   has negative curvature everywhere except at two
points (the Lagrange points).  A corollary
of this negativity is the uniqueness of the  $1/r^2$ figure eight,
 a complete symbolic dynamics for encoding   
the collision-free solutions, and the fact that collision solutions are
dense within the bound solutions. 
 
\vskip .3cm

{\bf 1. Introduction and Results.}

We   study  the  planar three-body problem
 with an  attractive $1/r^2$
potential.  According to the Lagrange-Jacobi
identity (eq. (3.7) below)
 every bounded  solution must have  
zero energy  and constant moment of interia $I$ , and conversely,
if an initial condition has zero energy and $\dot I (0) = 0$
then that solution is bounded.  
Setting  the moment of inertia $I$ equal to a constant  defines a three-sphere in
configuration space.  
Rotations act on this sphere according to the
Hopf flow so that the quotient of the three-sphere
by rotations is the two-sphere or {\it shape sphere}. See figure 1a. Points of this shape sphere 
represent oriented similarity classes of triangles.  Newton's equations,
for solutions with   $H = 0, \dot I = 0$, 
push down to  the shape sphere to yield a   
a family of second-order ODEs 
parameterized by the angular momentum.  These ODEs
have  singularities at the three points  representing
the three types of binary collisions.  
Upon deleting the collision points
we arrive  at dynamics on   the {\it pair of
pants}, -- the two-sphere minus three points.  
When the angular momentum is zero the resulting  dynamical 
system  is,   after a time
reparameterization, the geodesic flow for a certain \Ri metric on the
pair of pants.  This metric is the (reduced)  Jacobi-Maupertuis metric for energy
$0$.    
\proclaim Proposition 1. Endow
the pair of pants with the Jacobi-Maupertuis
metric (equations (3.9a,b) below).   
 Modulo rotations, translation, 
and scaling, the set of bounded 
zero-angular momentum solutions
for the $1/r^2$ potential three-body problem are in bijective correspondence 
with   geodesics
for this metric.   
The  metric is complete 
and its ends (the deleted neighborhoods of the three
binary  collisions) are asymptotic to Euclidean cylinders of positive radii. 

{\bf Proof:}  Section 3.

\vskip .3cm

See figure 1b for a depiction of the pair of pants.  

\eject

\vskip 0.1in

\epsfxsize=2.50in

\epsfbox{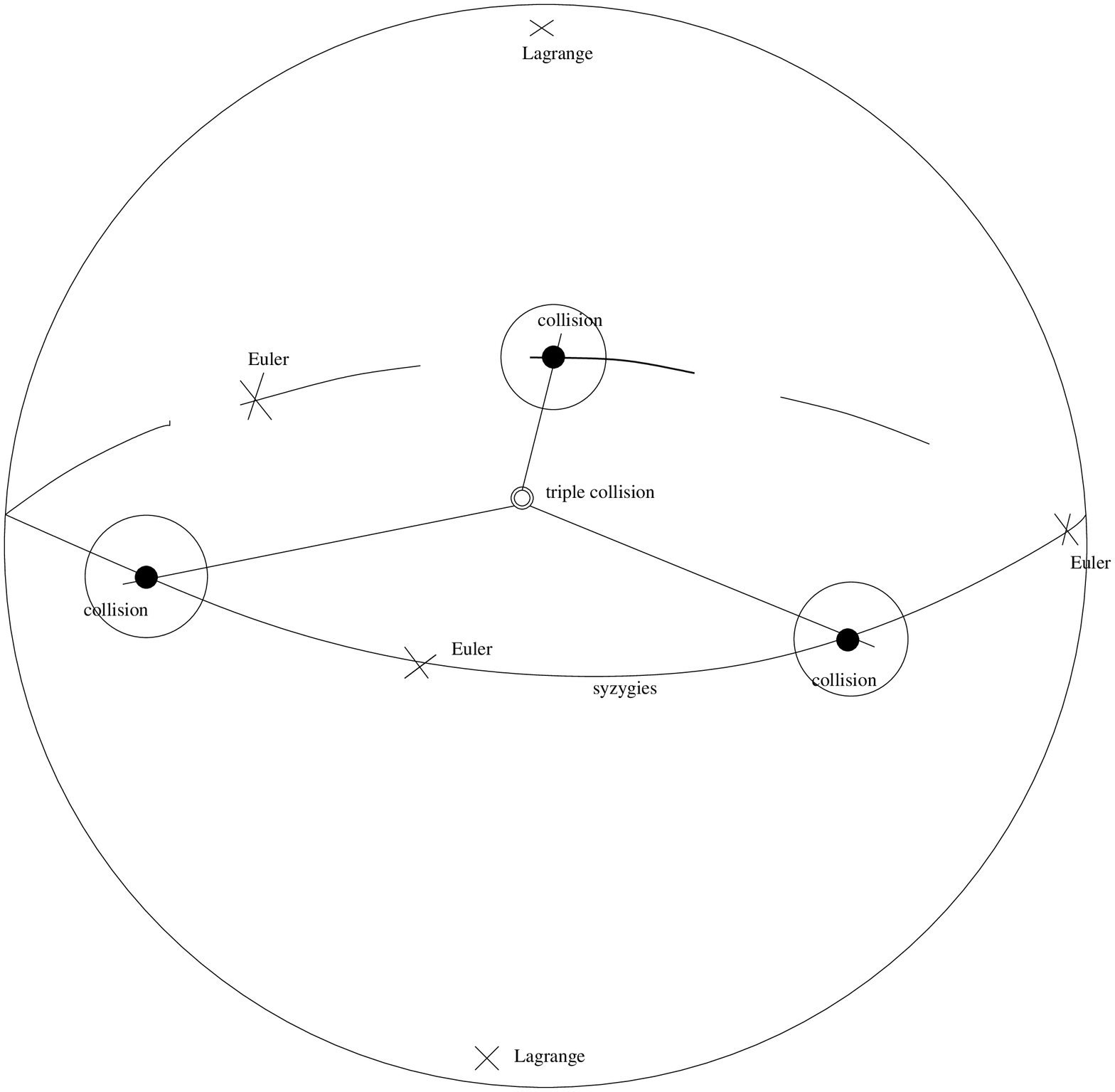}
\vskip 0.2in

{\bf Figure 1a.} The shape sphere 
\vskip .1in
\epsfxsize=3.00in

\epsfbox{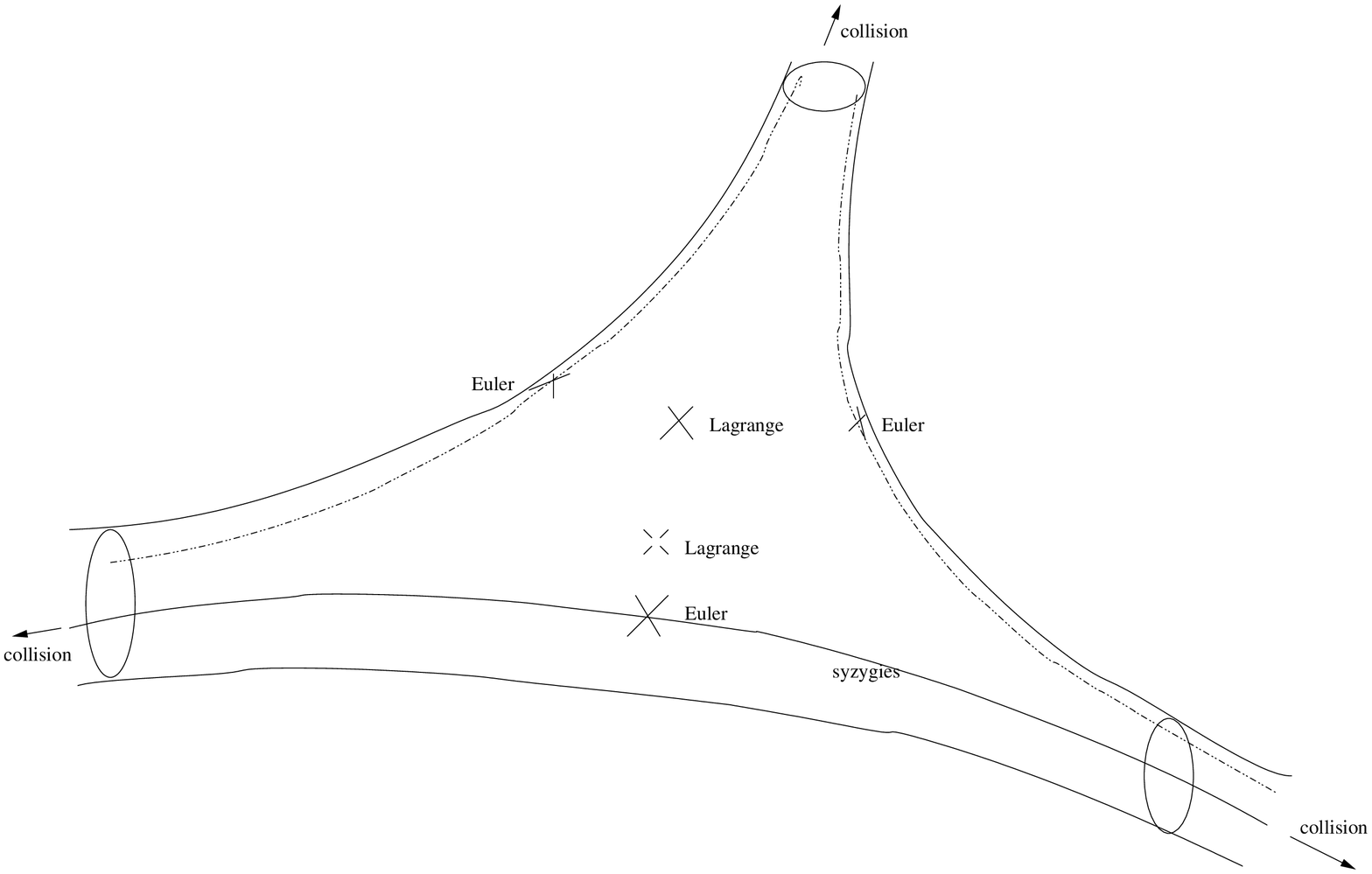}

\vskip 0.8in

{\bf Figure 1b.} The pair of pants

 The Jacobi-Maupertuis  metric depends parametrically on the   masses of the three bodies by way of
the potential (eq. 3.1).  Our main result is:
\proclaim Theorem 1.  If all three masses are
equal then  the
 Gaussian curvature for the Jacobi-Maupertuis metric on the pair of pants  is   negative 
everywhere  except at the two Lagrange points,  where it is zero.

{\bf Proof.}  Section 4. 

One might hope that  negativity   of 
the curvature persists for unequal masses.  It does not.  See section 7.
\vskip .44cm

{\bf 2. Motivation and  Dynamical Consequences.}  

{\bf 2.1. Periodic Orbits and their symbol sequences.}

This work began as  an attempt to give an analytic proof that  the
Newtonian ($1/r$ potential) figure eight solution of Moore-Chenciner-Montgomery
([Moore], [ChMont])
is unique.   
I began with the   easier
case of the $1/r^2$ eight.
 The figure eight  
is a periodic solution   which realizes a certain free homotopy class 
on the pair of pants.  Figure eights  exists for all $1/r^a$ potentials, $a > 0$
([CGMS] , [FerrTerr]).  For $a \ge 2$, 
 not only is the free homotopy class of the  eight realized,   but 
almost   every free homotopy class is realized by a solution. Combining these   facts
suggested the approach of this paper, and Theorem 1.

Our pair of pants metric from theorem 1 is neither compact, nor of negative curvature
everywhere.  But on a complete, noncompact surface
of negative curvature if a free homotopy class  has a geodesic representative,
then that representative is unique.  And  uniqueness
continues to hold if the curvature vanishes on a discrete set of points.
(This theorem is fairly well-known, and proved in  
in a more general context in section 6.4 below, and in particular
eq (6.4.3.).) We have proved
 
\proclaim Corollary.  For  the $1/r^2$ 
equal-mass zero-angular-momentum three-body problem,  if
a solution realizes a given  
free homotopy class on the pair of pants,
then that solution is unique   modulo rotation and scaling.
In particular the eight is unique modulo these
symmetries. 

In stating the corollary we begged the question of which classes 
are realized.  Every
class is realized  with the exception of  those classes which wind around a single end.
(See [MontN]).  
Gordon [1970] calls these  `bad' or unrealizable classes  `untied' 
while the complementary `good',  or realizable classes he  called `tied', being 
that they are 
`tied' to the collision singularities.  On the pants,
a bad class can be represented by  drawing a small circle,
or ``anklet''  around one
pants leg, and traversing it some number of times.
As this ``anklet''  is pushed   down towards the end of the leg
its length decreases.  As a result, any   minimizing sequence 
of curves realizing such a class ``falls off'' of the leg. (See theorem 3 below.)

We   follow [MontN] in using syzygies to describe the tied and untied classes.     
A {\it syzygy} 
is a collinear configuration of  the three bodies.
Syzygies come in three flavors,
1,2, and 3, depending on which mass is  between the other two.
(We exclude collisions.) 
The collinear configurations form the equator of the shape sphere.
(Figure 1a.) The three collisions lie on the equator so
that deleting them divides the  equator into three arcs,
again labelled 1,2,3  according to the   mass in the middle. 
   A  curve  on the 
shape sphere has  an 
 associated syzygy sequence:     list the  syzygies in order. 
(Assume that the   syzygy times  are discrete.)
The syzygy sequence of a motion of the three bodies
is obtained by projecting the motion onto the
shape sphere and writing out the syzygy sequence of 
the curve resulting  on the shape sphere.   
 
 Periodic  curves
give rise to periodic sequences.    For example, the class
in which  1 and 2 circle about each other for ever  while 3 remains far
away  has  syzygy sequence $\ldots 121212 \ldots$. 
(This is a ``bad'' class as we see later.)   We   subject  syzygy sequences to   the {\it  no stutterring rule}:
if  $ij$ are   consecutive letters
of the sequence, then   $i \ne j$.  The
reason for imposing this rule is that a stutter   can be homotoped
away. See figure 2. 

 A letter $j$ with a plus superscript,
as in $j^+$, denotes that syzygy $j$  occurs by crossing from the upper
to the lower hemisphere of the sphere.  A $j^-$
means that  
syzygy $j$  occurs by crossing from the lower
to the upper hemisphere of the sphere.   Pluses and minuses must alternate, since 
the  path will alternate between hemispheres.   In topological
terms an arc segment with two consecutive
pluses, such as   $k^-i^+ j^+$ lies entirely in
the upper hemisphere and   can be homotoped to $k^- j^+$.
In dynamical terms such an arc can never occur since at $i^+$
it would have to be tangent to the collinear subspace (the equator), but 
if a solution is tangent to the collinear  subspace at a point then
it lies completely within the collinear subspace.   

It follows from the above considerations, and the topology of the
pair of pants, that there is a one-to-one onto correspondence between  
 free homotopy classes and  periodic signed non-stuttering syzygy sequences.
From now on we will drop the signing indications --the $+, -$ superscripts -- for simplicity. 
(Given an unsigned sequence there are only two ways to decorate it with
signs.) 
The untied (unrealizable)  classes are precisely  those with syzygy sequence
  $...1212...$, $2323...$ or $...3131....$.
 They correspond to a curve  winding around a single end.  We prove in theorem 3 
that there is no bounded zero angular momentum solution which realizes
them.    
  Excluding these classes is equivalent to  insisting that
all three letters occur in the sequence.   Thus  the corollary asserts that  {\bf
every periodic non-stuttering syzygy sequence in which all three letters occur
is realized by a unique (up to symmetry) relative periodic solution.}

\vskip 0.2in

\epsfxsize=2.00in

\epsfbox{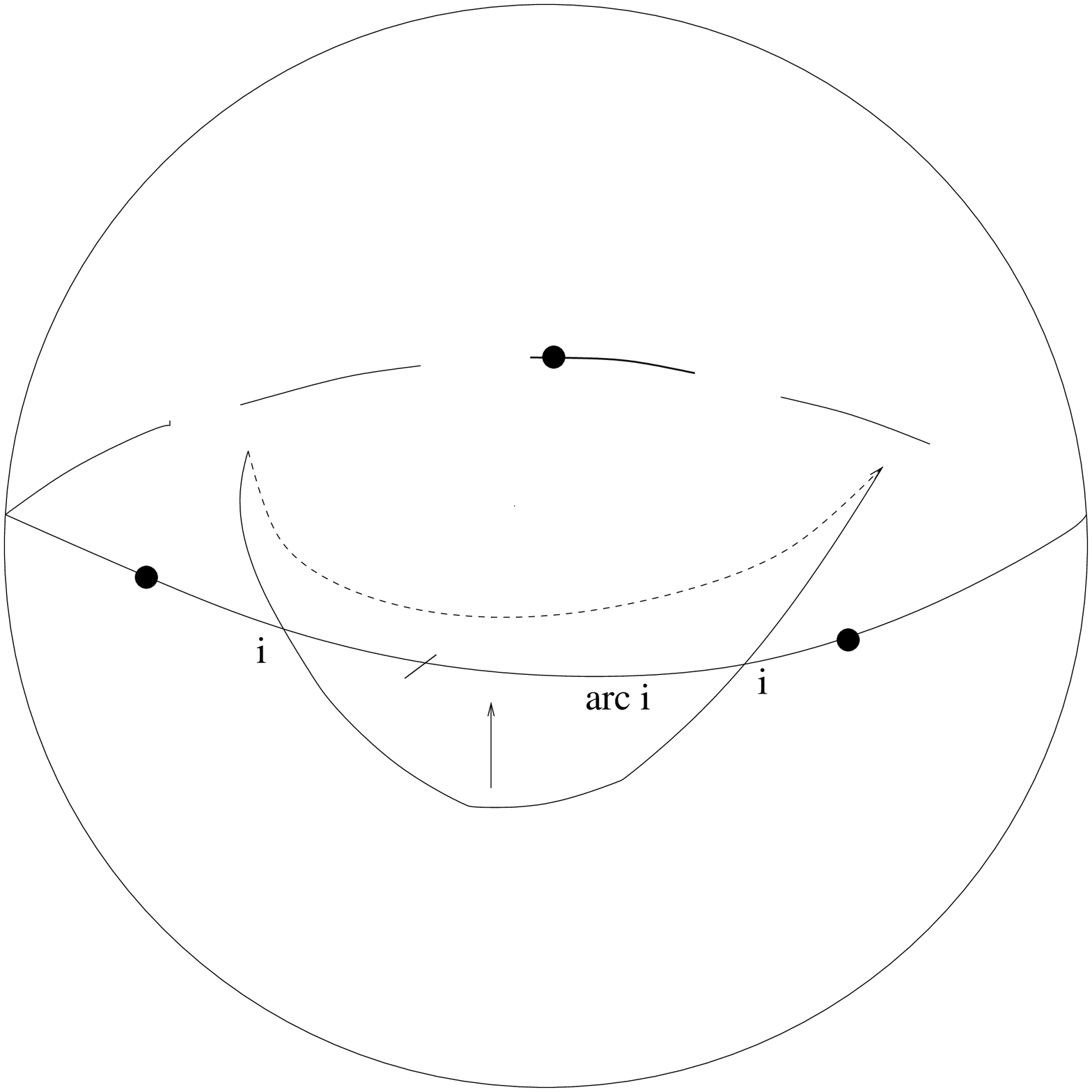}
\vskip 0.4in

{\bf Figure 2.} Homotoping away a stutter.

\vskip .4cm

{\bf 2.2. Symbolic Dynamics;  aperiodic syzygy  sequences.}
 
We  move  on to infinite aperiodic syzygy sequences. 
 In the rest of this subsection, `the problem'  means 
 the 
$1/r^2$ 
equal-mass zero-angular-momentum three-body problem, and
`solution' means a   solution to the problem, i.e. this differential equation.
Many of the ideas and results here are adaptations of those
pioneered in [Morse] and [Hadamard] XX Ref: use Had? .

\proclaim Theorem 2.  Every infinite nonstuttering syzygy sequence with
the exception of the untied classes $\ldots ijij \ldots$  is realized by
a solution. 

{\bf Proof.} Section 6.2.   The method   is the classical one 
[Morse] of approximation
by periodic solutions.   

\vskip .1cm

\proclaim Theorem 3.  If a  
syzygy  sequence ends (begins) with 
 $ijij\ldots$  then any bounded solution which realizes
this sequence must  end (begin)
in  the  $ij$ collision.  
The untied sequences $\ldots ijij \ldots$ are not
realized by   any solution. 

{\bf Proof.} Section 6.3.
 \vskip .1cm

It is perhaps worth remarking that if a solution suffers 
  collision 
then it does so in   finite Newtonian time, but 
infinite `Jacobi time''.  

Inspired by theorem 3, we call the sequences appearing there
``collision sequences''.  In more detail:

\proclaim Definition.  A 
  bi-infinite nonstuttering  syzygy sequence 
$s = \{s_j\}_{j = -\infty} ^{+\infty}$ 
is a forward  {\it collision sequence}
 if 
one of  its forward tails
$\{s_j\}_{j> N}$ contains only two letters.
Similarly, we have backward collision sequences.  A
collision sequence is one which is either forward
or backward collision sequence.  In the contrary case,  all three letters occur
in every tail, and the sequence is called {\it collision-free}.

\vskip .1cm 

Does every solution have a syzygy sequence?
If so, is this sequence unique?  In  [MontI]  I showed that
every bounded noncollinear zero-angular momentum solution to the Newtonian
three-body problem  suffers infinitely many
syzygies, provided the solution does not tend to triple collision,
{\bf and provided}  
binary collisions are counted as syzygies.  (See [Fuji] for another proof.)
  That proof  works verbatim for any  $1/r^a$
potential, $a > 0$, with the exception that
we must exclude   binary collisions. (They cannot
be regularized.)   Thus
every  bounded  {\bf collision-free} solution  has a syzygy sequence.  
That sequence must be nonstuttering in our equal mass case.
To prove that there is no stuttering,   use the fact that on a surface of negative
curvature any compact geodesic arc is the unique length minimizing curve among all
homotopic curves which share its endpoints. (See equation 6.3 and its derivation.)
Consequently, an application of  the method of reflection as exposed in
[ChM]  rids us of solution arcs representing stutters, i.e. solution arcs which
 hit  the same equatorial   arc
twice in a row.   These considerations allow us
to define a syzygy map from collision-free  sequences to 
infinite nonstuttering syzygy sequences. 

\proclaim  Theorem 4.   The
syzygy map from bounded solutions to   syzygy sequences
is a bijection between the set of   collision-free  solutions,   modulo
symmetry and time-translation,  and  the set of  bi-infinite nonstuttering collision-free syzygy sequences,
modulo shift.  

{\bf Proof.} Section 6.4.
 \vskip .1cm

Finally, we would like to know how much of phase space (the unit tangent
bundle of the pair of pants) is taken up
by the noncollision solutions.  Not much: 

\proclaim Theorem 5. Solutions tending to   binary collision are dense
within the space of all bounded solutions.  Thus the 
collision-free solutions have empty interior. 

{\bf Proof.} Section 6.5. 
 
\vskip .3cm 
{\bf Summary.} Putting  the theorems together gives a rather 
 complete symbolic dynamical
picture of the dyanmics of our problem -- the  zero-angular momentum
equal-mass $1/r^2$ three body problem restricted to the
bound orbits -- those with $I = const.$.  There are no linearly stable periodic orbits, by theorem 1.  
We will use the word ``bounded''
in the rest of this  paragraph to mean solutions  which tend to   to collision,
as these orbits are precisely the geodesics on the pair of pants which tend to infinity.
  
 Theorem 4 provides a  complete symbolic dynamics  
picture for the    bound  orbits : they are precisely the orbits none
of whose tails agree with the tied sequences $ \ldots ijij \ldots$.
The  closure of this set of orbits is the recurrent set.  The recurrence
set also   coincides  with  the closure of the
set of   periodic orbits.  There are unbounded orbits
on the frontier of this closure.   The situation is  similar to that of  
the recurrent set for the Kepler problem:  
  the space of periodic orbits is the recurrent set
and  contains
the unbounded parabolic orbits.    These unbounded  recurrent orbits   are ``just
barely unbound''  in that the collision condition $J_1 ^2 - (m_1 + m_2)  \le 0$  occuring
in  the appendix A, inequality (12A) is an equality on these orbits.    The complement of the
recurrent set   consists 
of orbits tending ``strongly'' to a binary collision. These
strongly colliding  orbits  form an open set. (Appendix
A.)  Finally,  the density result, theorem 6,  is an analogue of 
 what one would like to prove for the honest $1/r$ three-body problem:
that the set of solutions tending to infinity (via tight binary pairs)
is dense, for fixed energy and angular momentum.  M. Hermann
calls this density question  ``the oldest problem in dynamical systems
[Hermann] XX.

{\bf Loose ends.} 
There are some collision orbits
which we have  left out of symbol sequence considerations.
The collinear solutions are not accounted for.  (Collinear
solutions should either  have no syzygy sequence or a continuum
of `$i$'s as their sequence, depending on one's taste.)  There are exactly six collinear solutions, two for each of the
three collision arcs, the two being related by reversing orientation.    
There are also collision orbits which end in collision but without the bodies winding around  infinitely often.
They `head straight in' to infinity  down one of the pants legs and their 
corresponding syzygy sequences will  truncate in the forward direction, for
forward time collision.
The simplest of these truncated solutions are the isosceles solutions.  
Again,  there are six of these by the same
  counting as for colliner solutions.  The
isosceles solution  $r_{ij} = r_{ik}$   begins and ends  at the $jk$
collision, and has exactly one syzygy in between,  the Euler
point in which $i$ is at the midpoint of $j$ and $k$. 
Its syzygy sequence is the single letter 
 `$i$'. Interpolating between collinear and isosceles is a  
one-parameter family of  solutions whose  syzygy sequences truncate.    
These interpolating solutions are the $\lambda$-curves
of  eq. (3.13), the   curves of constant $\chi$, in
the $\lambda, \chi$ coordinate system  there.  
I do not know if the syzygy sequences of these  
solutions are finite, or one-sided infinite.

\vskip .3cm

{\bf Open questions.} {\bf 1.} Can two distinct collision orbits
share the same  syzygy sequence?  

{\bf 2.} Are there any solutions  besides isosceles which have a finite syzygy sequence?
If so, can any finite syzygy  sequence occur as the syzygy sequence of
some collision orbit?  

\vskip 1cm

\vskip 1cm

\vskip 1cm
 
 {\bf 3. Set-Up and Proof of Prop. 1} 

Write $x = (x_1, x_2, x_3) \in \R^6$ with $x_i \in \R^2$
for the positions of the three bodies,
and $r_{ij} = \|x_i - x_j\|$
for the distances between them.  The
potential is $-U$ where  
$$U = \Sigma m_i m_j/ r_{ij} ^2 \hskip 1cm (3.1).$$
The $m_i$ are the masses. 
Set
$$
\eqalign{ K  &= \Sigma m_i \| \dot x_i\|^2  \cr
&  = \langle \dot x , \dot x \rangle, \hskip 1cm (3.2) 
}$$
for twice the kinetic energy. The last
equality of (3.2) defines  the ``mass inner product'' on the
three-body configuration space $\R^6$.  The total energy 
$$H = K/2 - U \hskip 1cm (3.3)$$
 is constant along solutions.   The 
equations of motion, 
 $  \ddot x_i = -2 \Sigma_{i \ne j} m_j (x_i -x_j)/
r_{ij} ^4$, $i = 1,2,3$,  can be written
as the single vector equation
$$\ddot x = \nabla U \hskip 1cm (3.4)$$
where $\nabla U$ is defined using the
mass inner-product: 
$dU(x)(v) = \langle \nabla U(x), v \rangle$.

By the standard method of freshman physics, we can,
without loss of generality restrict our
considerations to motions for which
 $$\Sigma m_i x_i = 0 \hskip 1cm (3.5)$$
throughout.
This constraint defines a four-dimensional
real vector space which can be identified
with the two-dimensional complex space
$\C^2$ in such a way that counterclockwise rotation
of a triangle $(x_1, x_2, x_3)$
by $\theta$ radians  turns into scalar
multiplication of the corresponding complex vector by
$exp(i \theta)$. 
 Set  
$$ 
\eqalign{I & = \Sigma m_i m_j r_{ij}^2/ \Sigma m_i \cr
& = \langle x , x \rangle 
} \hskip 1cm (3.6)
$$
where the last equality is only true when 
the center of mass constraint (3.5) is in place.  
Using $\dot I = 2 \langle x ,\dot x \rangle$,
$\ddot I = 2 \langle \dot x, \dot x \rangle + 
2 \langle x , \ddot  x \rangle$ and
$\langle x, \nabla U (x) \rangle = - 2 U(x)$
(by $U$'s homogeneity) 
we obtain  the Lagrange-Jacobi identity: 
 $$\ddot I = 4 H \hskip 1cm (3.7), $$
valid along any solution. 
Thus 
$I(t) = const.$  along the solution
if and only if   $H = 0$  and $\dot I (0) = 0$
for that solution.

We will call a solution ``bounded'' if
the $r_{ij}$ are bounded as functions of time,
and  
do not simultaneously tend to zero, i.e. to triple collison.    
Now $I \to \infty$ if and only if one
 of   the $r_{ij}$ tend to infinity,
and   $I \to 0$ if all $r_{ij} \to 0$.
It follows from (3.7)  that
every bounded solution must satisfy $H = 0$,
$\dot I (0) = 0$, and $I(t) = const.$.
  
 The scaling symmetry
$x(t) \mapsto \lambda ^{-1/2}x (\lambda t)$
takes solutions to solutions, preserves zero energy,
and takes $I$ to $I/\lambda$.  Using this scaling, 
we may, without loss of generality,
assume that  $I = 1$ in studying   bounded solutions.
The set $I = 1$, $\Sigma m_i x_i = 0$
forms a three-sphere in the $\C^2$.  
{\bf We have reduced 
the study of the bounded solutions to the $1/r^2$ problem
to 
a second order dynamics on this three-sphere.}
It is well-known 
([Arn] or [AbMar]) that a 
dynamics on a constant energy surface 
$H  = E$
level set is equivalent to   
geodesic flow for the Jacobi-Maupertuis metric 
$( E + U) ds^2$
where $ds^2$ is the kinetic energy metric.
In our case, $E = 0$, and 
we restrict the kinetic energy to the sphere $I = 1$.  Consequently, 
the study of bounded solutions is
equivalent to the study of geodesics on
the three-sphere under the metric 
$ds^2_J =  U ds^2$
conformal to the standard metric $ds^2$ on
that three-sphere.   
 
To obtain a metric
on the shape sphere, we quotient
by rotations.   We review the discussion
in [MontN], [ChMont], or [MontR] on
this metric. See especially the appendix
of [MontR] for explicit computations 
and derivations.  The group of rigid rotations acts  on the three-body
configuration space   
according to scalar multiplication on $\C^2$ by unit modulus complex scalars.
Restricting ourselves to   
the three-sphere $S^3:= \{ I =1 \}   \subset \C^2$ and forming 
the quotient by this rotational action yields the
famous Hopf fibration 
$$S^3 \to S^2 = S^3/ S^1 \hskip 1cm (3.8) .$$
The quotient two-sphere is
the  shape sphere ([ChM], [MontR] esp. the appendix).
Its points represent oriented similarity classes
of triangles. 
Both the  dynamics and  the  Jacobi metric $U ds^2$
on the three-sphere   descend under the projection (3.8)
to the shape sphere  
once we fix the value of the total angular momentum.   
The total angular momentum of a solution is zero if and only if
that solution is orthogonal to the rotational orbits,
i.e. orthogonal to the fibers of (3.8).  The projection of 
such a zero-angular momentum  
 solution under (3.8)
is  a geodesic for the quotient metric.
(See [Hermann], lemma 4.1. His situation is  more general
than ours. Our circle  bundle (3.8) is replaced
by a general Riemannian submersion.) 
We can write the Jacobi-Maupertuis metric on the
shape sphere as  
$$ds^2_J = U ds_{shape}^2 \hskip 1cm (3.9a)$$
 where $ds_{shape}^2$ is
the kinetic-energy induced metric on shape sphere,
and $U$ is the (negative) potential (3.1) restricted
to $I =1$ and then viewed as a fucntion on the shape-sphere
(possible because of its rotation invariance.
The shape sphere metric is the   round metric on a sphere of radius $1/2$:
$$ds^2 _{shape} = ({1 \over 2})^2 [d \phi ^2 + \cos^2 (\phi) d
\theta ^2) ] \hskip 1cm (3.9b)$$
where
 $\phi$ is the colatitude --
the angle from the equator, and $\theta$
is labels longitudinal circles on the sphere.
\vskip .4cm
\vskip .4cm
{\bf Proof of Proposition 1.} The 
discussion of the last two paragraphs shows
that,  modulo rotations, translation, 
and scaling, the set of bounded 
zero-angular momentum solutions
for the negative potential (3.1) are in bijective correspondence 
with   geodesics
for the Jacobi-Maupertuis metric on the shape sphere
minus the three binary collision.   
 Under this correspondence the
geodesic flow for the metric corresponds,
after a time reparameterization, to the
the flow defined by Newton's equations. 
It remains to verify the claims about the completeness
and  that the ends asymptote to cylinders.
  
\vskip .4cm

{\bf Completeness.}

Let $\rho = \rho_{ij}$ be the spherical
distance from  the $ij$ collision point $C_{ij}$,
as measured in the spherical metric  $ds^2_{shape}$.
Then as $\rho \to 0$ we will show that 
$$U = {C^2 \over \rho^2} + O (1)  \hskip 1cm
(3.10a) $$
for some positive constant $C^2$, while 
$$ds^2 _{shape} = d \rho^2 + (\rho^2 + O(\rho^4)) d \chi^2
 \hskip 1cm (3.10b)$$
where  $\chi$ is  the  angular coordinate based at $\Sigma_{ij}$ 
so that $(\rho, \chi)$ are geometric polar
coordinates.   It follows that the Jacobi metric has
the expansion:
$$ ds^2_J =   {C^2 \over \rho^2} (d \rho^2 + \rho^2 d \chi^2)) + O(1) \hskip 1cm
(3.10c).$$ It follows
that if we  approach the collision $\rho = 0$ along any curve then the
length of that curve diverges at least as fast as the integral
of $C\sqrt{d \rho^2/ \rho^2} = Cd\rho/\rho$, that is,  it  
diverges logarithmically as $C|\log(\rho)|$ as $\rho \to 0$.
Consequently any  curve tending towards ``infinity''
i.e. to one of the binary collisions,   has infinite length,
which proves completeness.

To establish (3.10a), it suffices to establish
$$r_{ij} = {1 \over {\sqrt{\mu_{ij}}}}  \sin(\rho_{ij})
\hskip 1cm (3.11)$$
where $\mu_{ij} = m_i m_j/(m_i + m_j )$ is the reduced mass.
The other two distances $r_{ik}, r_{jk}$ are bounded
away from zero as $r_{ij} \to 0$, due to 
the constraint $I =1$
(see (3.6)).  Then (3.10a) follows from (3.1)
  and
the Taylor expansion of $\sin(\rho)$.  The constant
$C$ in (3.10a) is $m_i m_j / \sqrt{\mu_{ij}}$.

 To establish (3.11) we work in 
the full three dimensional shape space which is the space whose points are
oriented congruence classes of planar triangles.  The full shape space is 
isometric to the   cone over the shape sphere and consequently
distances $d$ in the full shape space can be obtained 
from spherical distances together with knowldge of
the distance $R = \sqrt{I}$  from the cone point. 
Write $d_{ij}$ for the distance in the full shape space between
an arbitrary point and   the $ij$ binary collision {\it ray}.
Then we have  
$$r_{ij}  = {1 \over {\sqrt{\mu_{ij}}}} d_{ij} \hskip 1cm (3.12a) $$
and
$$d_{ij} =  R \sin(\rho_{ij}) \hskip 1cm (3.12b)  .$$
Upon setting $R= 1$, (3.11) follows immediately.
Equation (3.12a,b) can be 
 found in section 4, equations (4.3.15a,b) of [MontN]. However, note that  there
is a typo in eq 4.3.15a. The $\mu_{ij}$ in that equation must be replaced by
$\sqrt{\mu_{ij}}$.)   
 
To get (3.10b)  use the fact that the shape sphere is isometric
to the sphere of radius $1/2$ and that the metric on such a sphere is given by 
$$ds^2_{shape} = d \rho^2 + [(1/2) \sin(2\rho))]^2 d\chi^2 \hskip 1cm (3.13).$$
in spherical-polar coordinates. Then use the Taylor expansion of
$(1/2) \sin(2 \rho)$. 

\vskip .4cm 

{\bf Asymptotes to Cylinders.}  We 
use a more precise version of the expansion (3.10c).
Set 
$$d\lambda = - \sqrt{U} d \rho \hskip 1cm (3.14).$$
Integrating (3.14) defines   a  function 
$\lambda = \lambda(\rho, \chi)$ such thta 
 $\lambda \to \infty$
as  the collision $\rho = 0$ is approached. 
From (3.10a) we have  $d\lambda = - C d \rho / \rho  + O(1)$
from which it follows that 
$$\rho = e ^{ -C \lambda} + o(\rho).$$
The  $(\lambda, \chi)$ are coordinates for the end $\rho = 0$,
and in these coordinates
$$d\bar s^2 = d \lambda^2 + f(\lambda, \chi)^2 d \chi^2
\hskip 1cm (3.15)  $$
where, from (3.13) and (3.9a) we have  
$$f^2  = ({1 \over 2} \sin (\rho) )^2 U.$$
Now $({1 \over 2} \sin  ( 2\rho) )^2 =  \sin^2 \rho \cos^2 \rho$
so that from (3.11)  
$ ({1 \over 2} \sin (2 \rho) )^2 =  \mu_{ij}
r_{ij}^2 \cos^2 (\rho )$
 and
$$f^2 =  \mu_{ij}  \cos^2 (\rho)\{ m_i m_j  + m_i m_k {{r_{ij}^2} \over
{r_{ik}^2}} 
 +m_j m_k {{r_{ij}^2} \over {r_{jk}^2}} \} \hskip 1cm (3.16)$$
where $ijk$ is a permutation of $123$. 
As we approach the collision $\lambda = \infty$
we have $r_{ij} \to 0$ while $r_{ik}, r_{jk}$ remain bounded
since we are constrained to $I =1$.  Thus
$$\lim_{\lambda \to \infty} f = \sqrt{\mu_{ij}  m_i m_j}: = K_{ij} > 0 \hskip 1cm (3.17).$$ 
Summarizing:  
$$d   s^2 _J  = d \lambda^2 + (K_{ij} ^2 + O( e^{-2 C  \lambda})) d \chi^2 
\hskip 1cm (3.18)$$ 
which says the metric asymptotes to a Euclidean cylinder of
radius
$K_{ij}$ as we approach the $ij$ end.  

 QED 
\vskip .1cm
{\bf Remark.}  It follows from equations (7.12, .13)
that   the Gaussian curvature near the end is
negative.  But for any metric
of the form (3.15) this curvature is equal to $ -{1 \over f}{{\partial ^2
f}
\over {\partial \lambda ^2}}$.  Thus, for 
fixed $\chi$,   the 
function $f(\lambda, \chi)$ is a  strictly convex
of $\lambda$, for all $\lambda$ from some point on, 
and from this point on,   
$f(\lambda, \chi)$  monotonically decreases to $K_{ij}$.
\vskip .4cm
 
\vskip 1cm 

{\bf 4.  Curvature. Proof of theorem 1.}

We proceed to the proof of our main result, Theorem 1,
the negativity of the Gaussian curvature   when the masses are equal.
The computation   proceeds 
through a series of lemmas.  The first  is standard and we will
not provide the proof. 

\proclaim Lemma 4.1.  Let a  surface 
be endowed with   conformally related  metrics
$ds^2$ and $d \bar s^2  = U ds^2$.
Then their curvatures $K, \bar K$ 
are related by
$$\bar K = U^{-1} (K - {1 \over 2} \Delta \log(U))$$
where the Laplacian $\Delta$ is with respect to the
$ds^2$ metric.

The curvature $K$ of the standard shape metric  $ds^2 = ds^2 _{shape}$  of (3.9b) is 
$K = 4$. According to lemma 4.1
$$\bar K = 4 - {1 \over 2} \Delta (\log(U) \hskip 1cm (4.1)$$
is the desired curvature, the curvature
 of  the  Jacobi-Maupertuis metric (3.9a) 
of proposition 1 and Theorem 1.
A routine computation yields
$$\Delta (\log(U)) = {{ U \Delta U - \|\nabla U \|^2} \over
{U^2}} \hskip 1cm (4.2) .$$
Here, and throughout this section,  $U$ is considered  
  as a function on the shape sphere,  
$\Delta U$ is its Laplacian with respect to the standard
shape space metric metric $ds^2$  
   and  
$\| \nabla U \|^2$  is the squared length of its gradient with
respect to the same metric. 

A key to the subsequent computations
is to use the squared length coordinates as in [AlbCh] 
$$s_k = r_{ij}^2,  ijk \hbox{ a permutation of } 123,  \hskip 1cm (4.3)$$
rather than the lengths $r_{ij}$ themselves. 
Write 
$$U_{2n} = \Sigma 1/r_{ij}^{2n} = \Sigma 1/s_k^n \hskip 1cm (4.4)$$
so that $U = U_2$.
\vskip .2cm

\proclaim Lemma 2.
  $$\Delta U = 8U_4 \hskip 1cm (4.5) .$$

{\bf Proof. } Section 5.2.

\vskip .2cm

\proclaim Lemma 3.
$$\| \nabla U \|^2 = 4 S \hskip 1cm (4.6a)$$
where
$$S = 2U_6 -  U_4 -{3/2} \Sigma^{\prime} 1/s_i ^2 s_j ^2 +  2 \Sigma^{\prime} 1/s_i  s_j ^2
- \Sigma^{\prime} 1/s_i  s_j  
\hskip 1cm (4.6a)$$
and where ``$ \Sigma^{\prime}$ '' means to sum over all indices 
$i, j$ with $i \ne j$. (For example  $ \Sigma^{\prime} s_i s_j = 2s_1 s_2 + 2 s_2 s_3 + 2 s_3 s_1$,
 twice the second symmetric polynomial in the $s_i$.)

{\bf Proof. } Section 5.3.

\vskip .2cm

{\bf Proof of the Negativity of the curvature.}
Combining the equations (4.1), (4.2), (4.5) and (4.6a) 
we find 
$$ -  \bar K U^3 =  4 U U_4 -  4 U ^2  -  2 S \hskip 1cm (4.7).$$
Expand out the first two terms on the right hand side:
$$\eqalign{ U U_4 & = \Sigma 1/s_i \Sigma 1/s_j^2 \cr
&= \Sigma 1/s_i ^3 + \Sigma^{\prime} 1/s_i s_j ^2 \cr
& = U_6 + \Sigma^{\prime} 1/s_i s_j ^2 .
} \hskip 1cm (4.8) 
$$
while
$$\eqalign{ U ^2 & = \Sigma 1/s_i \Sigma 1/s_j \cr
&= \Sigma 1/s_i ^2 + \Sigma^{\prime} 1/s_i s_j  \cr
& = U_4 + \Sigma^{\prime} 1/s_i s_j  .
} \hskip 1cm (4.9) 
$$ 
Plugging  (4.8), (4.9) and equation (4.6b)  
 back into equation  (4.7) yields :  
$$ 
-  \bar K U^3   = -2 U_4  - 2 \Sigma^{\prime} 1/s_i  s_j  + 3  \Sigma^{\prime} 1/s_i ^2  s_j ^2 
\hskip 1cm (4.10a).$$
Use the fact that $U_4 + \Sigma^{\prime} 1/s_i  s_j = (\Sigma 1/s_i)^2 = U^2$ 
to rewrite the right hand side of (4.10), and divide the resulting equation
 by two in order to  obtain
$$ 
- \bar K U^3   = 3  \Sigma^{\prime} 1/s_i ^2  s_j ^2  - 2 U^2  \hskip 1cm (4.10b)
$$
Consequently, $\bar K \le  0$ if and only if  
$$3(\Sigma^{\prime} 1/s_i ^2 s_j ^2) \ge 2 ( \Sigma 1/s_i)^2 \hskip 1cm (4.11)$$
To prove (4.11), multiply both sides of it by $s_1 ^2 s_2 ^3 s_3 ^2$
thus arriving at $6 (\Sigma s_i ^2) \ge 2 \sigma_2 ^2$ or
$$3 (\Sigma s_i ^2) \ge  \sigma_2 ^2 \hskip 1cm (4.12)$$
where $\sigma_2 = s_1 s_2 + s_2 s_3 + s_3 s_1$
is the second elementary symmetric polynomial in the $s_i$.
(The coefficient $6$ arose because for each pair $ij$ there
are two terms in the sum $\Sigma^{\prime}$.  See the parenthetical
remark in  lemma 3.)
To prove (4.12), remember that we are restricting ourselves
to the sphere $I = 1$ and that $I = \Sigma s_i /3$.  Thus we 
can homogenize the equation by using that $3 =  (\Sigma s_i )^2 /3$ 
on the sphere.  So, the desired inequality now reads:
$${{(\Sigma s_i )^2} \over 3}(\Sigma s_i ^2) \ge  \sigma_2 ^2.$$
The  two inequalities:
$$\Sigma s_i ^2 \ge \sigma_2 \hskip 1cm  \hskip 1cm (4.13A)$$
and 
$${{(\Sigma s_i )^2} \over 3} \ge \sigma_2 \hskip 1cm (4.13B) $$
are classical, with 
equality  in {\bf either case} if and only
if all the $s_i$ are equal.   Here are the proofs. 
 Inequality (4.13A) follows simply upon 
rearranging  the inequality 
$(s_1 - s_2)^2 + (s_2 - s_3)^2 + (s_3 - s_1)^2 \ge 0$.
 Inequality (4.13B) is a special case of a general
inequality among the elementary symmetric polynomials evaluated at 
positive arguments $s_i$.  See for example the 
Encyclopaedia [Math],  
App. A, Table 8, inequality (4).  Alternatively, expand
out $(\Sigma s_i)^2 = \Sigma s_i^2 + 2 \sigma_2$
and use  inequality (4.13A).
Multiplying the two inequalities
yields the desired inequality (4.12),
with equality if and only if we are at the Lagrange points
$s_1 =s_2 = s_3$ of the shape sphere.
Since (4.12) is equivalent to the curvature inequality, 
Theorem 1 is proved, modulo the proofs of lemmas 2 and 3
which follow in the next section.   QED

\vskip 1cm

{\bf 5. Proofs of lemmas 2 and 3.}

{\bf 5.1. Notation.}

To compute $\Delta U$ and
$\|\nabla U \|^2$ we must express
the squared distances $s_k = r_{ij}^2$ of (4.3) 
in terms of the  spherical coordinates $(\phi, \theta)$ of (3.9b). 
In [MontI] I prove that upon restriction to the sphere $I = 1$ 
$$s_k = 1 - \cos(\phi) \gamma_k (\theta) \hskip 1cm (5.1.1b)$$
 where 
$$\gamma_1 (\theta) = \cos(\theta), \gamma_2 (\theta) =
\cos(\theta + 2 \pi/3), \gamma_3 (\theta) = \cos(\theta + 4
\pi/3) \hskip 1cm (5.1.2).$$ The special angles $\theta = 0, 2 \pi/3, 4 \pi/3$
mark the locations of the three binary collision
on the equator $\phi = 0$ of collinear triangles. 
Later on we will use the 
fact that the three planar vectors
$(\gamma_k, \gamma_k^{\prime})$, $k =1,2,3$ form the vertices of an equilateral
triangle inscribed in the unit circle.  
\def\dphi{\partial_{\phi}}
\def\dtheta{\partial_{\theta}}
{\bf Here and throughout we   write  $\gamma_i ^{\prime}$
for the derivative $\dtheta \gamma_i$
of $\gamma_i$ with respect to $\theta$.}

\vskip .2cm  

{\bf 5.2. Proof of Lemma 2.} 
Write $c = \cos(\phi)$,
$\dphi$ for the partial derivative with respect
to $\phi$ and $\dtheta$ for the partial derivative
with respect to $\theta$.  Then 
$$\Delta U = {4 \over c} \dphi (c \dphi U) + {4 \over c^2} (\dtheta
^2 U). \hskip 1cm (5.2.1)$$
And
$$\dphi (1/s_i) = -s\gamma_i /s_i ^2 \hskip 1cm (5.2.2a)$$
$$\dtheta(1/s_i) = c\gamma_i ^{\prime}/s_i ^2 \hskip 1cm (5.2.2b). $$
Since  
$U = \Sigma 1/s_i$
we have
$$
\eqalign{
\dphi(c \dphi U)  & = \dphi c \Sigma (-s\gamma_i)/s_i ^2
\cr
&= \dphi(-cs \Sigma \gamma_i)/s_i ^2  \cr
     &= (-c^2 + s^2)\Sigma \gamma_i /s_i ^2 
 + 2c s ^2 \Sigma \gamma_i ^2/ s_i^3}.
 $$
Thus
$$
{1 \over c}
\dphi(c \dphi U) = ((s^2 - c^2)/c)\Sigma \gamma_i /s_i ^2 
 + 2 s^2 \Sigma \gamma_i ^2/ s_i^3 .\hskip 1cm
\hskip 1cm (5.2.3)$$
And
$$\eqalign{{1 \over c^2} \dtheta \dtheta U 
& = {1 \over c^2} \dtheta \Sigma c \gamma_i ^{\prime}/s_i ^2
\cr
& = {1 \over c^2} \Sigma c [\gamma_i ^{\prime \prime}/s_i ^2
+ 2 (c \gamma_i )^2/ s_i ^3
\cr
& = {1 \over c} \Sigma \gamma_i ^{\prime \prime}/s_i ^2
+ 2 \Sigma (\gamma_i ^{\prime})^2/s_i ^3
}
\hskip 1cm (5.2.4)$$
Now $\gamma_i^{\prime \prime} = - \gamma_i$
so that
$${1 \over c^2} \dtheta \dtheta U =
-{1 \over c} \Sigma \gamma_i /s_i ^2 
 + 2 \Sigma (\gamma_i ^{\prime})^2/s_i ^3 \hskip 1cm (5.2.5)
$$

Adding 4 times (5.2.4) to 4 times (5.2.5) we get
$$\Delta U = 
4 ((s^2 - c^2 - 1)/c ) \Sigma \gamma_i /s_i ^2  
+  8 s^2 \Sigma \gamma_i ^2/ s_i^3   + 
8 \Sigma (\gamma_i ^{\prime})^2/s_i ^3 \hskip 1cm (5.2.6) $$
Now
$s^2 -c^2 - 1 = -2c^2$ so that
$$\Delta U = -8 \Sigma c \gamma_i /s_i ^2  
+ 8 \Sigma [s^2 \gamma_i ^2 + (\gamma_i ^{\prime})^2]/s_i ^3
$$
Recalling that  $\gamma_i ^2 + (\gamma_i ^{\prime})^2 =1$
(the vectors $(\gamma_i,  \gamma_i ^{\prime})$
define an equilateral triangle inscribed in the unit circle)
we see that we can replace  $(\gamma_i^{\prime})^2$ by $1 - \gamma_i ^2$
in order to  obtain
$$\eqalign{
s^2 \gamma_i ^2 + (\gamma_i ^{\prime})^2
 & = (s^2 -1) \gamma_i ^2  +1 \cr
& = - c^2 \gamma_i ^2 +1   \cr
& = -(1-s_i)^2 +1  \cr
& = 2s_i -s_i^2 \cr
& = s_i (2 -s_i)
}
\hskip 1cm (5.2.7) $$
Then 
$[s^2 \gamma_i ^2 + (\gamma_i ^{\prime})^2]/s_i ^3
= (2 -s_i)/s_i^2 = (c \gamma_i +1 )/s_i^2$
where I used $1 -s_i = c \gamma_i$
It follows that
$$
\Delta U = -8  \Sigma c \gamma_i /s_i ^2  
+ 8 \Sigma c \gamma_i /s_i ^2 +  8 \Sigma (1/s_i ^2)
= 8U_4
$$
as claimed.
 
\vskip .4cm
{\bf 5.3. Proof of Lemma 3.}

We have  
$$\eqalign{
dU  &= \dphi U d \phi + \dtheta U d \theta \cr
   & = ( \Sigma (-s \gamma_i /s_i ^2  )d\phi + ( \Sigma (c \gamma_i ^{\prime})/s_i^2) ) d \theta
} \hskip 1cm (5.3.1) 
$$
Now $\| \nabla U \|^2 = \|d U\|^2$.  The length squared of
the covector $dU$ is computed relative to the metric `$g^{ij}$ '
induced on covectors, which, from (3.9) is 
given at the  point $(\phi, \theta) $
by  $\|a d\phi + b d \theta \|^2 = 4( a ^2 + {1 \over c^2} b^2)$,
with $c = \cos(\phi)$.
It follows that
$$
\eqalign{
\|\nabla U \|^2  & =  4(\Sigma s \gamma_i/s_i^2)^2 +  4 (\Sigma \gamma_i ^{\prime}/s_i ^2 )^2 
\cr
& =  4 \big( \Sigma s^2 \gamma_i ^2 /s_i^4 + \Sigma^{\prime} s \gamma_i s \gamma_j / s_i^2 s_j^2 
+ \Sigma \gamma_i ^{\prime ^2} /s_i^4 + \Sigma^{\prime} \gamma_i ^{\prime} \gamma_j ^{\prime} / s_i ^2 s_j ^2 
\big)
\cr
& = 4 ( \Sigma (s^2 \gamma_i ^2 + (\gamma_i ^{\prime 2} ) /s_i ^4 +
\Sigma^{\prime} (s^2 \gamma_i  \gamma_j + \gamma_i ^{\prime} \gamma_j ^{\prime})/ s_i^2 s_j^2 ).  
} \hskip 1cm (5.3.2) 
$$
Simplify the numerator
in the first summand of the last equation
by  using (5.2.7). 
 We will simplify the numerator
of the second summand  by using  an analogous identity for
$s^2 \gamma_i  \gamma_j + \gamma_i ^{\prime} \gamma_j ^{\prime}$, $i \ne j$.
Indeed, since the vectors $(\gamma_i, \gamma_i ^{\prime})$
form the vertices of an equilateral triangle inscribed within
the unit circle,we have that 
$\gamma_i \gamma_j + \gamma_i ^{\prime} \gamma_j ^{\prime} = -1/2$
for $i \ne j$, since  $-1/2 = \cos(2\pi/3)$ is the cosine
of the central angle defined by any two vertices of an equilateral triangle. Thus 
$$\eqalign{
s^2 \gamma_i  \gamma_j + \gamma_i ^{\prime} \gamma_j ^{\prime}
&= \gamma_i  \gamma_j + \gamma_i ^{\prime} \gamma_j ^{\prime} - c^2 \gamma_i \gamma_j
\cr
&=  -1/2 - (1- s_i)(1-s_j) \cr
&= -3/2 + s_i + s_j - s_i s_j \hskip 1cm  
}. \hskip 1cm (5.3.3)
$$
Plugging (5.3.2) and (5.3.3) into (5.3.1) we get
$$
\eqalign{
\|\nabla U \|^2  & =  4 ( 2 \Sigma 1 /s_i ^3 - \Sigma 1/s_i ^2 
-{3 \over 2} \Sigma^{\prime} 1/ s_i ^2 s_j^2 + 2 \Sigma^{\prime} 1/s_i s_j ^2  - \Sigma^{\prime} 1/s_i s_j )  
\cr
&= 4 S}
$$
as claimed. 
QED

\vskip .4cm

{\bf 6.  Dynamical Consequences.} 

{\bf 6.1}

Let $P$ denote the pair of pants.  
Using the spherical shape metric (3.9a),
construct three disjoint circles with centers at
the three binary collision points.  Delete the open discs
bounded by these circles to obtain a   compact region $R \subset P$ 
having for its   boundary  the three disjoint circles.   Refer back to  figure 1.

 \proclaim Lemma 6.1. Any  arc in $P$ whose (finite) syzygy sequence 
contains all  three letters   must cut through  $R$.

{\bf Proof.} Let $c$ be such a `123'  arc.
If one of $c$'s syzygies is within $R$, we are done.
Otherwise, all three syzygies lie within the three excised discs.
 But  all three syzygies cannot be   in the same disc,
since each disc contains exactly  two syzygy types. Thus
$c$ must   must travel  from one disc to the other, and in 
 so doing it cross into $R$.  
 
QED.

\vskip .3cm

{\bf 6.2. Proof of Theorem 2.} Let $s$ be 
syzygy sequence containing all three letters.
Approximate $s$ by a sequence $w^N$, $N = 1,2, 3$ 
of periodic sequences as follows.   
Truncate   $s$ to form the finite even length
subword 
$w_N = s_{-N +1} s_{-N +1} \ldots s_N$.   
Turn this subword into a periodic sequence
$ \ldots w_N w_N w_N \ldots$ by repeating it
in blocks.  If the resulting word
has stutters at the join,  shift the ``window''
we used to form $w^N$ so as to form the word 
$w_{N, j} = s_{-N + j +1} s_{-N + j} \ldots s_{N + j}$
along with its  corresponding  periodic word.  
We can always find a $j$ so that the resulting periodic word,
call it $w^N$,  is non-stuttering.
The sequence 
$w_N$ 
contain all three letters $123$ for all $N$
sufficiently large, since $s$ itself contains all three letters.  
Theorem 2 implies that  the $w_N$, for $N$ large,  
are  represented by a unique geodesic
$\gamma_N$. By lemma 6.1 the $\gamma_N$  must cut through $R$.  
Shifting time (and thus shifting the sequence) if necessary, we may assume that
$\gamma_N (0) \in R$.  Since $R$ is compact, 
so is  the unit tangent bundle to $P$   
(relative to the Jacobi-Maupertuis metric) over $R$.   
The pairs
$(\gamma_N (0), \dot \gamma_N (0)$ lie
in this compact space, so we  
we can find a subsequence of them   which 
  converge to some initial condition $(q,v)$.
The geodesic with initial condition $(q,v)$
realizes the   infinite sequence
$s$. This establishes the existence of a  solution 
realizing the   syzygy sequence $s$.

QED

\vskip .2 cm

FIGURE, SECTION AND EQ RELABELLING NEEDED  BELOW   XX

{\bf 6.3. Collision sequences. Proof of theorem 3.} 

Let $s$ be a collision sequence.  
We may  assume, 
without loss of generality,  that it is a forward
collision sequence, and that the two letters in its  
forward tail are  $1$ and $2$.  The first part of the theorem
asserts that  any solution  which realizes $s$
must satisfy $r_{12} \to 0$ as
the Jacobi time  $t \to \infty$.
(In Newtonian time  the  collision  occurs in finite time.
As the     two bodies get closer,
the third body  effects them less and less.
A straightforward  analysis of the  two-body $1/r^2$
problem shows that 
$\lim \inf r_{12} = 0$ if and only if
$\lim r_{12} =0$.  The perturbation
of the third body does not affect this assertion.
Thus in order to prove the solution realizing $s$ suffers
collision it suffices to show that  
$$\lim \inf r_{12} = 0 \hskip 1cm (6.3.1) .$$
along  the solution.

Our proof of (6.3.1)  relies on the fact that on a simply connected complete surface
of  non-negative curvature any compact geodesic arc is the
 unique minimizer between
its endpoints.   On  a complete non-simply connected surface such
as the pair of pants   this implies
that if we have a geodesic arc, then there is no shorter curve
which share  endpoints with that arc and which is homotopic to it
through endpoint-fixing homotopies.    

We argue by contradiction.  Suppose that some
solution $\gamma \subset P$ realizes the collision sequence
$s$ but satisfies 
$\lim \inf r_{12} = \delta > 0$.
  We will construct a comparison curve  $c$ which
has the same endpoints as a (long) arc of $\gamma$,  is homotopic
to this arc 
through end-point fixing homotopies, 
but which is shorter than $\gamma$.  The  arc will
be one whose syzygy sequence is $1212 \ldots 12$ with
$N$ repeats of $12$, and $N$ large.  
 See figure 6 for the picture of this arc
and the  shorter comparison curve $c$.  The existence
of $c$ contradicts the minimality of $\gamma$ described in  the previous 
paragraph.

To  construct $c$ we will use 
the cylindrical  coordinates
$(\lambda, \chi)$ 
of  (3.14a,b), (3.15) associated
to the  $12$-collision end.  We have    
$$\lambda (\rho, \chi) = \int_{\rho_0} ^{\rho} \sqrt{U(s,
\chi)} ds$$
where $(\rho, \chi)$ are spherical-polar coordinates
centered at the collision. 
 The $(\lambda, \chi)$ coordinates are valid 
on the entire sphere minus the `$3$' equatorial arc
and the collision points.
 The coordinate $\lambda$ satifies 
  $\lambda  \to \infty$  as collision
is approached.  The Jacobi metric in these coordinates is 
$$ds^2 = d \lambda^2 + f(\lambda, \chi)^2 d \chi^2
\hskip 1cm (6.3.2)  $$
where  
$$f^2 = {1 \over 2}   \cos^2 (\rho) \{ 1   +
 {{s_{12}} \over
{s_{13}}} + {{s_{12}} \over
{s_{23}}} \} \hskip 1cm  $$
The curvature of any metric of the form (6.3.2) 
is  
$ -{1 \over f}{{\partial ^2 f} \over
{\partial \ \lambda ^2}}$.  Since this
curvature is negative (theorem 2) and since 
$\lim_{\lambda \to \infty} f(\lambda, \chi) = 1/\sqrt{2}$
(eq. 3.18) 
we have that   for each fixed
$\chi$ the function    
$f(\chi, \lambda)$ decreases monotonically to its infimum $1 / \sqrt{2}$
as $\lambda \to \infty$.  
It follows that 
$F(\lambda) = \min_{\chi} f(\lambda, \chi)$
also   decreases monotonically to $1/\sqrt{2}$.

Any geodesic arc $\gamma$ on the pair of pants which
realizes the syzygy sequence $12$  
cannot cross either isosceles circle
$r_{12} = r_{23}$ or
$r_{12} = r_{13}$.  This follows from the minimality 
property of the arc, discussed above, and   the reflection principle 
(as in [ChM] and the proof of no stuttering in section 2, between
theorems 3 and 4).
Reflections about the isosceles circles are isometries of the Jacobi metric,
so that any segment of $\gamma$  which crosses, then crosses back, can be reflected,
so as to form a new arc with the same endpoints as $\gamma$, and the same homotopy
type, contradicting uniqueness.   Thus, without loss of generality,
we may assume that our long subarc of $\gamma$ lies entirely in
the union of the regions $r_{12} < r_{13}$ and
$r_{12} < r_{23}$.   In particular,
the $(\lambda, \chi)$ coordinates are valid all along our
arc.  

To construct the desired comparison curve $c$ we will  
use the fact that $F$ is monotone
decreasing in $\lambda$  and that the number $N$ of $12$ crossings of  a subarc of $\gamma$
can be taken arbitrarily
large.      Since $\inf r_{12} > 0$, by assumption, we have
that $\Lambda = \sup \lambda < \infty$ along our curve.
Let $\epsilon > 0$ small be given. Choose  points $\gamma(s_N), \gamma (t_N)$
$s_N < t_N$, along the arc 
for which $\lambda (t_n), \lambda (s_n)  >  \Lambda - \epsilon$ 
and such that  the arc $\gamma[s_n, t_n]$ in between realizes
the syzygy sequence $12 \ldots 12$ with $N$ copies of $12$.
Our comparison curve $c$ will go ``straight in'' to collision
until some point with $\lambda_* > \Lambda$, to be determined momentarily,
winds around the $12$ collision point in the same sense as $\gamma$ for the same number of
syzygies at this fixed value $\lambda = \lambda_*$, and then return headed ``straight out''
from collision to the point $\gamma(t_n)$.  See figure 6.  ``Straight in'' and ``straight out''
means that $\chi$ is fixed, and only $\lambda$ varies. During
its `winding around'  journey, $\lambda = \lambda_*$ is
fixed and $\chi$ varies, starting at  $\chi =
\chi(\gamma(s_n))$, increasing so that $c$ 
suffers the   syzygy sequence $1212 \ldots 12$ (N
times), and then stopping   at $\chi = \chi(\gamma(t_n)$ in time for the ``straight out''
return segment.

\vskip 0.1in

\epsfxsize=3in

\epsfbox{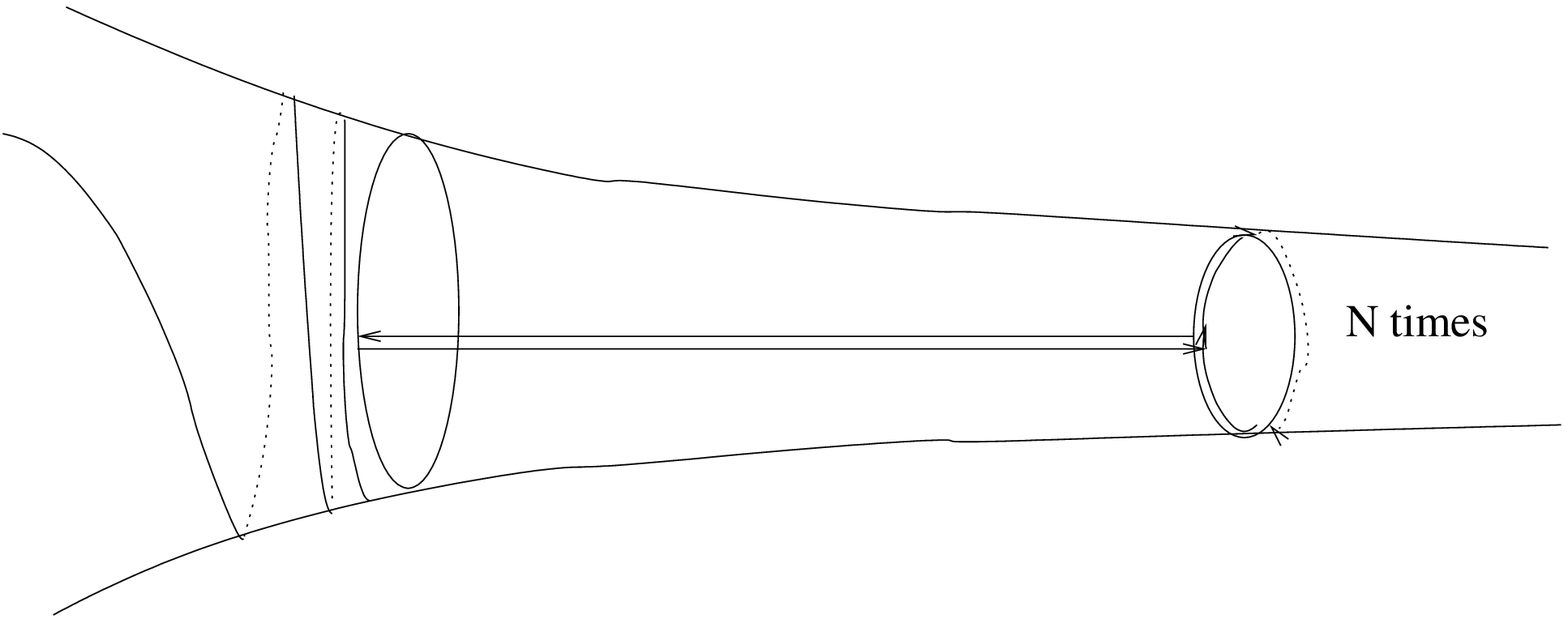}
\vskip 0.2in

{\bf Figure 6.} Orbit surgery to shorten length of
a collision sequence path

FFXX LABEL FIG as per next to last par of this sec.

\vskip .1in
 
\vskip 0.8in

By construction, the curve  $c$ shares endpoints with our arc of $\gamma$,
and is homotopic to it.
It remains to  show that $c$ is shorter  than our arc of $\gamma$.
By the monotonicity and the limiting properties of $f$ we can 
choose $\lambda_* > \Lambda$ so that 
$$max_{\chi}f(\lambda_*, \chi) < F(\Lambda):  = min_{\chi} f(\Lambda, \chi).$$
Choose $N$ so large that
$${{2(\lambda_* - \Lambda) + \epsilon + \pi F(\Lambda) } \over N} < \pi (F(\Lambda) -
max_{\chi}f(\lambda_*,
\chi) ).$$  It follows that 
$$ \pi N  max_{\chi} f(\lambda_*, \chi) + \pi max_{\chi} f(\lambda_*, \chi)   + 2(\lambda_*
- \Lambda) + 2 \epsilon <  N\pi (F(\Lambda)$$
Now consider the cartoon in figure XX or the
construction of $c$.   The `in and out' arcs of
$c$ have length less than $2(\lambda_* - \Lambda) + 2 \epsilon$.
The `around arc' between the $N$ syzyygies has length  less than $\pi N  max_{\chi}
f(\lambda_*,
\chi)$.  (Refer again to figure 6.) The additional $\pi max_{\chi} f(\lambda_*, \chi)$
accounts for the fact that $\chi(s_n)$ and $\chi(t_n)$ need
not be equal.      Thus the  left hand
side of the inequality is greater than the length of $c$.
A similar but simpler analysis shows that the right hand side is smaller than
the length of our arc of $\gamma$.   Thus $c$ is shorter than the arc of $\gamma$,
completing  the proof of first assertion of  theorem 3. 

The same analysis shows that the infinite sequence $\ldots 1212 \ldots$
of all $12$'s is never realized.  For such a realization must
be a local minimizer, and the above orbit surgery shows we can always
decrease the length of a path by making it closer to collision.

QED

\vskip .3cm 

{\bf 6.4. Proof of theorem 4.} 

To establish the   uniqueness of the realizing solutions, we 
work on the
universal cover $D$ of the pair of pants $P$. Topologically, $D$ is  
 the Poincare disc and the  
fundamental group  $\Gamma = \pi_1 (P)$ (the free group on two letters) 
acts on $D$ as a Fuchsian group.    See figure 3,
  and also the
book Indra's Pearls [Mumf].  
In  figure 3b   we have drawn in a fundamental domain $P_0$
for $P$ and some of its images   under $\Gamma$. The
boundary of $P_0$ consists of four circular arcs which are lines 
lines relative to the Poincare metric on $D$.  These bounding 
arcs are labelled $1_+, 1_-$ and $2_+, 2_-$. To form $P$
out of $P_0$ glue 
arcs $1_+$ and $1_-$ to form syzygy arc $1$, and glue
arcs $2_+$ and $2_-$ to form syzygy arc $2$.
Syzygy arc $3$ is internal
to the fundamental domain.  In the figure we have dropped the
$+, -$ subscripts on the arcs. 

\eject

\vskip 0.1in

\epsfxsize=3.00in

\epsfbox{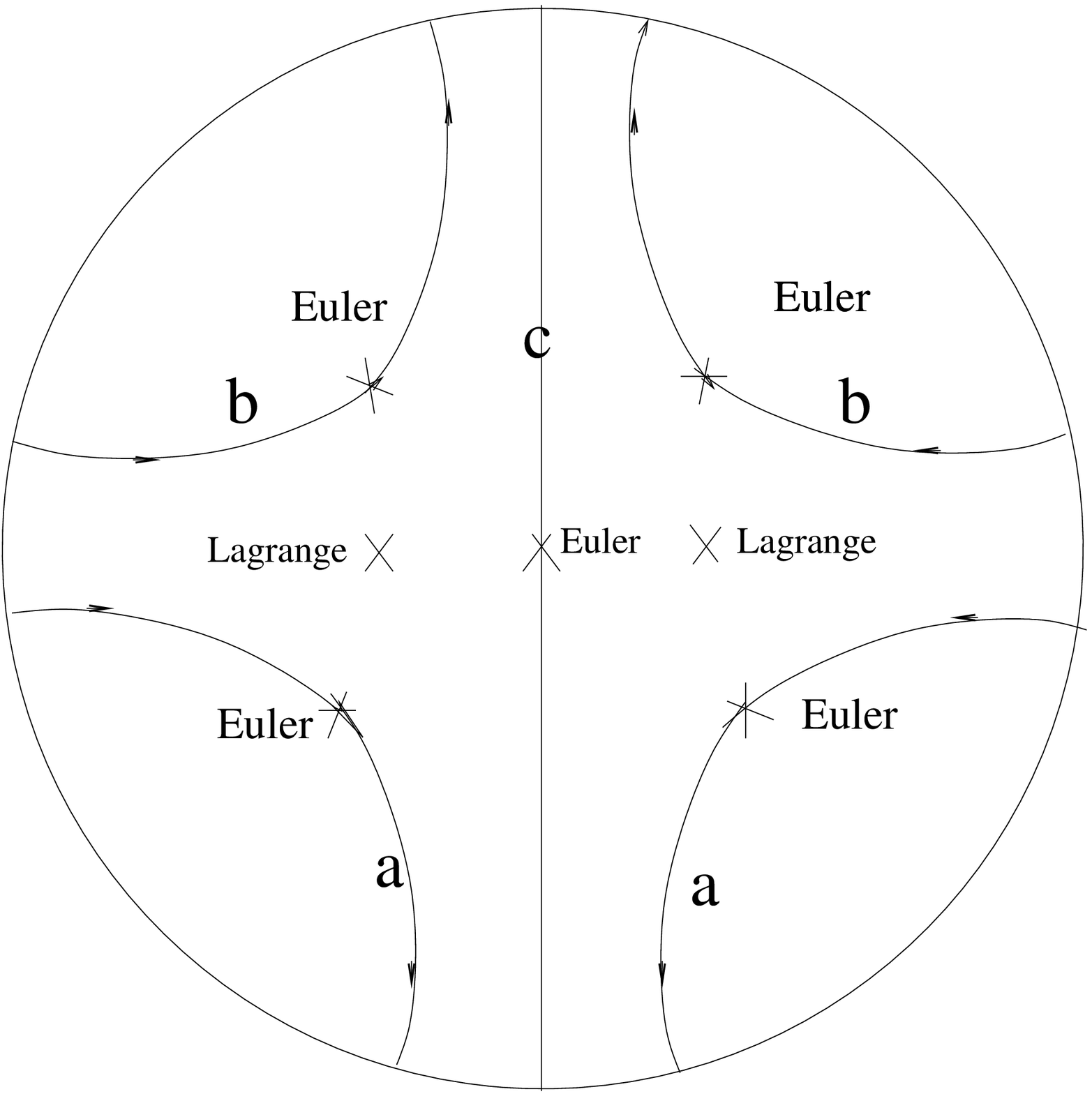}
\vskip 0.1in

{\bf Figure 3a.} The fundamental domain
\vskip .2in
\epsfxsize=3.00in

\epsfbox{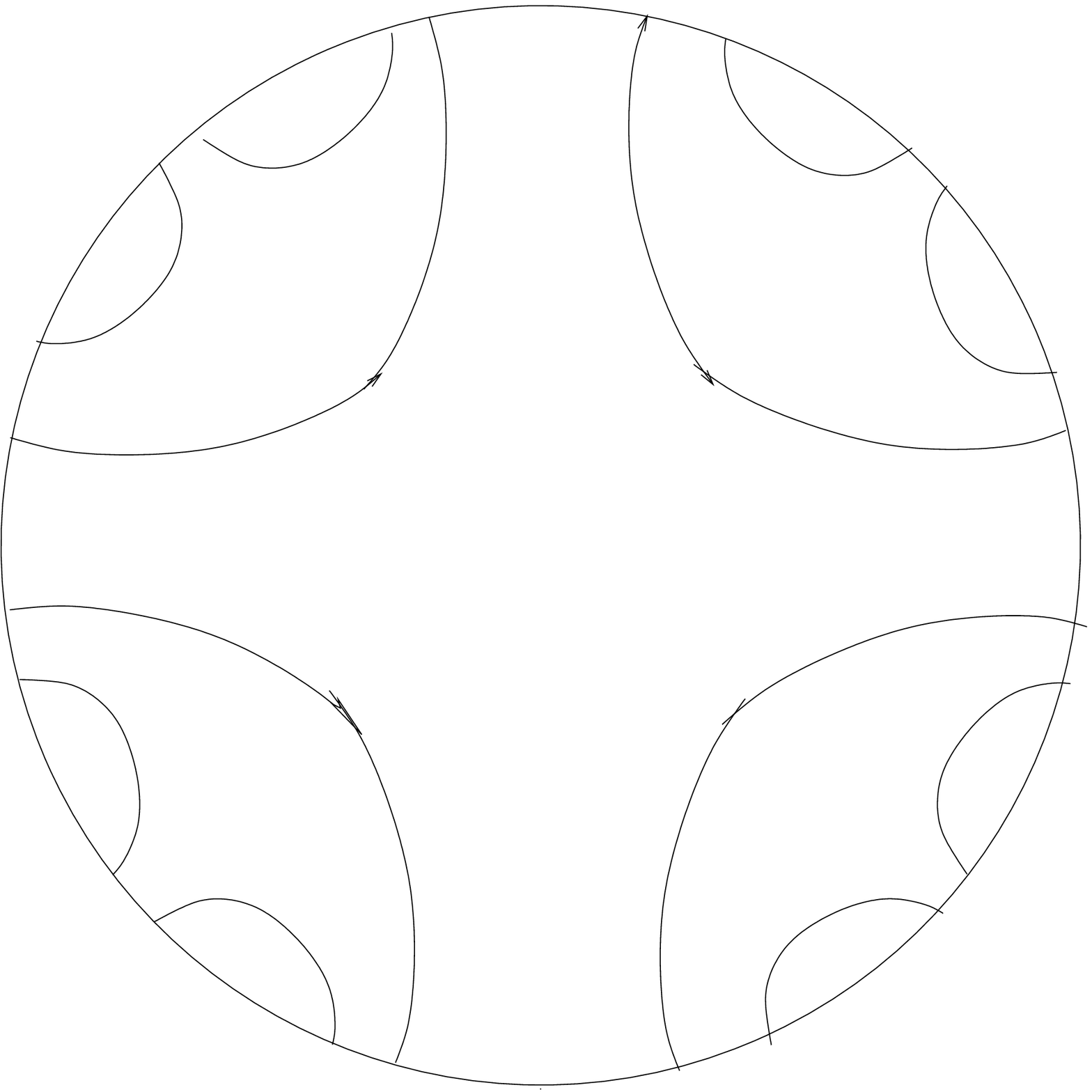}

{\bf Figure 3b.} Some of the tiles.

\vskip 0.4in

Write $$\pi: D \to P$$
for the covering map, so that the fibers of $\pi$
are copies of $\Gamma$.   We will use the hyperbolic, or   constant 
negative curvature metric on $D$ 
in order to understand the fundamental group $\Gamma$ and its action on $D$.
(With   respect to this   metric  $\pi$ is not a local isometry.)  
$\Gamma$ acts on $D$ with respect to hyperbolic
isometries, i.e. M\"obius transformations. 
$\Gamma$ is freely generated by two elements $a$ and $b$,
one of which, say $a$,  interchanges $1_+$ and $1_-$, and 
the other of which, $b$, 
interchanges $2_+$ and $2_-$. These elements act on $D$
as Mobius transformations.  (Viewed
as acting on the Riemann sphere, they  interchange exteriors with interiors of 
their respective circles. )
 In order to form $P$ out of $P_0$ glue arc $1^+$ to $1_-$ 
by $a$  
and glue  $2_+$ to $2_-$ by $b$.  
The projection under $\pi$ of these boundary arcs form the arcs $1,2$ of the equator
of $P$.
The third arc $3$ is internal to $P_0$, and separates it into
two halves, the northern (+) hemisphere, and southern (-) hemisphere.
      
 The images of $\gamma P_0$ of the fundamental domain 
$P_0$ under elements of $\gamma \in \Gamma$ tile all of $D$. 
Each of the four  boundary arcs of such a tile $\gamma P$  
 is the image under $\gamma$ of a unique
boundary arc of $P_0$ and we continue to  label the tile's boundary arcs  
by the corresponding $P_0$ lables  $1_+, 1_-, 2_+$, or  $2_-$. 
Two  tiles intersect, if at all,   along 
a common boundary arc.  This common
arc must be   a   $j_-$ arc of one tile
and a $j_+$ of the other,   $j =1, 2$. 
   
Suppose  now  that 
 two geodesics   realize the same symbol sequence $s$.
Denote by $\gamma,   c$   the lifts of these geodesics to 
  the universal cover $D$.  After  translating these  curves  by   elements of $\Gamma$ 
we may suppose that  both  begin in the  reference fundamental domain
$P_0$.   
I claim that the syzygy sequence $s$ uniquely specificies
a sequence of contiguous tiles  $\ldots P_{-2}P_{-1} P_0 P_1 P_2 \ldots $
through which   $c$ and $\gamma$ must pass.  
To see this fact, we first note that  each sequence of
three contiguous tiles $P_{i-1} P_i P_{i+1}$
represents either two or three letters of a signed syzygy sequence. See figure 4. 
The sequence is obtained by drawing a curve which crosses from $P_{i-1}$
through $P_i$ and into $P_{i+1}$ in the `most direct'' way.
The curve must enter into $P_i$ across one
of its bounding arcs $1_+, 1_-, 2_+,
2_-$. The choice of $P_{i-1}$ uniquely specifies
which arc.  It must leave   
across another such arc and the choice of  $P_{i+1}$
uniquely specifies this exit arc.  Along the way it must either
cross $3$ or not.   If no   internal syzygy with arc 3 occurs then the  sequence has two letters
and no `$3$', and otherwise the sequence
does contain the  letter $3$ as the middle letter.    
    Consequently,
both $\gamma, c$ pass through an identical
list of tiling domains, as claimed.  

Each tiling domain 
$P_j$ has within
it an inverse image $R_j = \pi^{-1} (R) \cap P_j$  of our  compact domain $R$.  
By lemma 6.1 both $c$ and $\gamma$
must have the property that for infinitely many $j$ we
have that both  $\gamma$ and $c$ 
lie in $R_j$. 

\vskip 0.2in

\epsfxsize=3.00in

\epsfbox{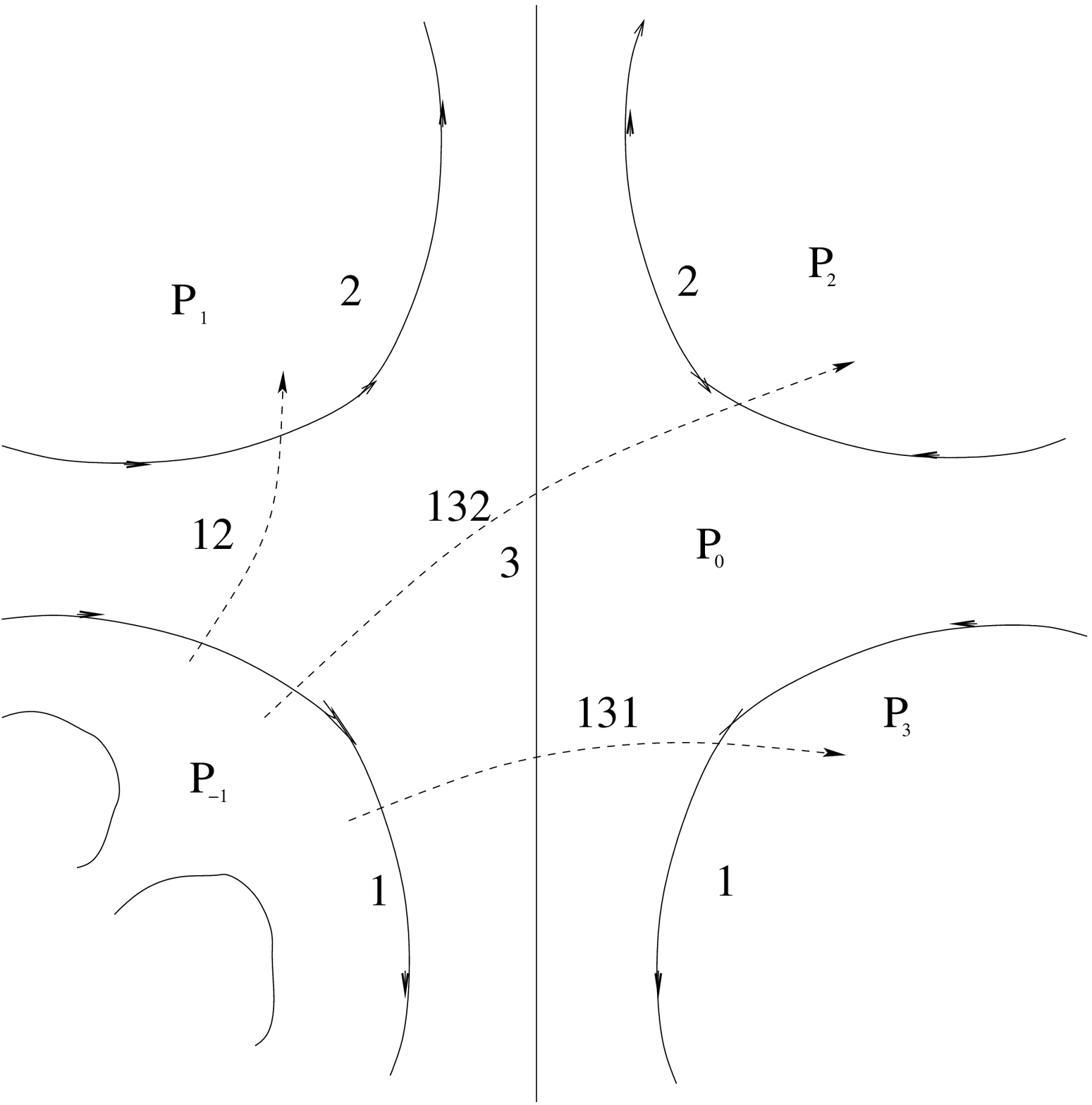}
\vskip 0.4in

{\bf Figure 4.} Syzygies and Tile Crossings.
\vskip .2in

 Now  lift the Jacobi-Maupertuis metric
from $P$ up to $D$, using the
projection $\pi: D \to P$, thus arriving at a complete
non-negatively curved $\Gamma$-invariant metric
on $D$ for which our two curves
$\gamma$ and $c$ are geodesics.
The curvature of this metric is   zero only
at the discrete set of points $\pi^{-1}(L_{\pm})$,
where $L_{\pm}$ are the Lagrange points of $P$.
And the map $\pi$ is  a local isometry for this metric.   
Write 
$$h(t) = dist(\gamma, c(t)) \hskip 1cm (6.4.1)$$
for the distance between the variable point $c(t)$ on
the curve $c$ 
and the entire geodesic $\gamma$.  
See figure 8.  Now
each $R_j$ has finite diameter $\delta$ because
$R$ is compact, and the   two curves pass
through the $R_j$'s infinitely often, indeed every time 
the letters  $123$ occurs contiguously in $s$. It follows that 
  
$$\lim \inf_{t \to +\infty} h(t) \le \delta 
\hbox{ and } 
\lim \inf_{t \to -\infty} h(t) \le \delta. \hskip 1cm (6.4.2)$$ 

We will now show that inequality (6.4.2) is impossible unless  
 the two geodesics are in fact the same, in which case $h(t) =0$ everywhere.
We use the formula
$$d^2 h /dt^2 = -sin(A(t)) \int_{d_t} K ds \hskip 1cm (6.4.3)$$
proved in the following paragraph.
In this formula, 
$d_t $ is  the geodesic realizing the distance
$h(t)$.  It  has one endpoint at the point $c(t)$ on $c$
and the other endpoint on $\gamma$ which it intersects
perpindicularly.  See figure 5. The angle $A(t)$ is the
angle of intersection
between the geodesics $c$ and
$d_t$ at  $c(t)$.  See figure 5. 
This angle satisfies   $0 < A(t) < \pi$, so that $\sin(A(t) > 0$  
The negativity of   $K$ 
(except at a discrete  point set) implies that 
$h(t)$ is strictly convex: $d^2 h /dt^2 > 0$.  
   But any strictly convex function defined
on the real line  tends to infinity  in
one direction or the other. This  contradicts   (6.4.2).
Our two geodesics must be the same.

\vskip 0.2in

\epsfxsize=3.00in

\epsfbox{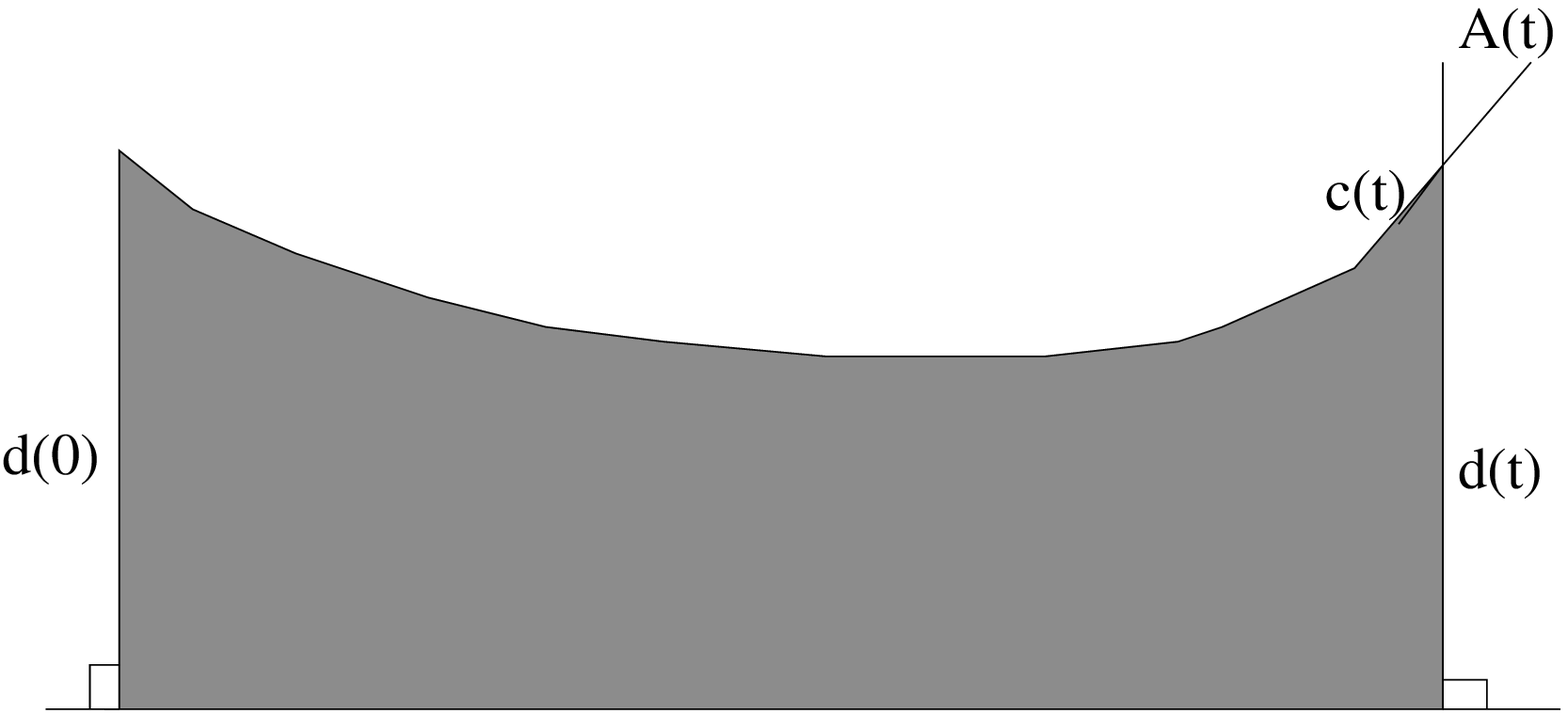}
\vskip 0.4in

{\bf Figure 5.} Variation of Distance
\vskip .2in
 
{\bf Derivation of (6.4.3). }

The first variation of arclength implies
 $$dh/dt = \cos(A(t)) \hskip 1cm (6.4.4) .$$
 Let
$M(t) \subset D$ denote the quadrilateral whose  edges 
consist of the geodesic arcs $d_0, d_t$
together with  the arcs of $\gamma, c$ which connect
$d_0$ to  $d_t$.  According to  the Gauss-Bonnet theorem,
for any such geodesic quadrilateral $Q$ we have
$$ 2 \pi - (\Sigma \hbox{ interior angles }) = - \int \int_Q  K dA$$
In  the case of $M(t)$ the interior angles are $\pi/2, \pi/2, \pi - A(0)$ and $A(t)$.  
See figure 5 again.  
Thus  
$$A(0)  - A(t)= - \int \int_{M(t)} K(t) dA \hskip 1cm (6.4.5) $$
Differentiating (6.4.4) with respect to $t$ yields  
$$dA/dt = \int_{d_t} K ds. \hskip 1cm (6.4.6) $$
Now differentiate (6.4.4) with respect to $t$, using (6.4.6)  to obtain (6.4.3).
 
\vskip 1cm

\magnification=\magstep1 \hoffset .65pt
\def\R{I\!\! R}

\def\C{C \!\! \! \!  I\, }

\def\RP{I\!\! R I\!\! P}
\def\Ri{Riemannian }

{\bf Proof of theorem 5.}  Let $c$ be a bound orbit and
$s = \{ s_j \}_{j = - \infty} ^{\infty}$ its syzygy sequence.
Approximate $s$ by a family $w^n$, $n = 1,2, 3, \ldots$  
of forward collision sequences by replacing the tails $s_j, j > n$ of $s$
by that of the $12$ collision sequence. Collapse two letters if neccessary when
stutters appear at the `join' $j = n$ of   the replacement. 
The backward tails of
the $w^n$ contain
all three letters because $s$ is collision-free, so  we can apply theorem 2 to realize the $w^n$ by
solutions $\gamma_n$ (not necessarily unique).  By lemma 6.1
the $\gamma_n$ all pass through $R$, so, shifting time if necessary,
we have that the tangent vectors   $v_n = (\gamma_n (0), \dot \gamma_n (0))$ 
are unit vectors with position $\gamma_n (0)$ in   $R$.  
We now argue as in the proof of theorem 3.  By compactness
of the  set of unit
tangent vectors over $R$, we can form a convergent
subsequence of the $v_n$, which we  relabel as  $v_n$, so that $v_n \to v$.
Let $\gamma$ be the curve with initial condition $v$.  
The curves $\gamma_n$  converge, over compact
sets,  to  $\gamma$, so that the syzygy sequences $w^n$ must converge to
the syzygy sequence of $\gamma$.  Thus $\gamma$ and $c$ share the same
syzygy sequence.  By   theorem 4,  $\gamma = c$.  Consequently the
unboounded curves $\gamma_n$ converge to our initial bound curve $c$.

QED

\vskip .3cm

{\bf 8. Curvature for  unequal masses.} 

In this section is to prove that Theorem 1 is
a special case:  for
most mass distributions the curvature takes on both signs.  
Take the masses to be general positive numbers
$m_i$.  Writing  $$p_i = m_j m_k, s_i = r_{jk}^2$$  for
$ijk$ a   permutation of $123$ we have
$$U = \Sigma p_i /s_i, \hskip .3cm  I = {1 \over M} \Sigma p_i s_i,
\hbox{  
where } M = \Sigma m_i.$$
The Jacobi metric is obtained by multiplying the
shape metric $ds^2_m$ on the shape sphere $I =1$  
by $U$ and restricting it   to the sphere $I =1$.
(The subscript `m' indicates dependence on the masses.)
It is conceptually and computationally more straightforward to identify
the shape sphere  as the space of rays in shape space
and to make $U$ a function on the space of rays by
making it  homogeneous of degree  0  by multiplying it by $I$. Thus, setting
$$\tilde U = I U , $$  our Jacobi metric is
$$ds^2_J =  \tilde U ds^2_m.$$

In order to compute, we will use the  coordinates  $\phi, \theta, R$  of   [Mont2] 
for  shape space. 
(See, in particular, the notation and computations of 
section 7 there.)
Write $I_1$ for the moment of inertial when
all masses are equal to one:
$$I_1 = \Sigma r_{ij}^2/ 3. \hskip 1cm (7.1) $$
The   $\phi, \theta$ are spherical coordinates
for the   $I_1$ shape sphere, while the radial coordinate  
$R = \sqrt{I}$ for the mass-dependent   $I$ of (3.6). 
 With respect to these coordinates we have that  
$$s_k  =   I_1 (1 - \gamma_k (\theta) \cos(\phi)  \hskip 1cm (7.2).$$
as before
Note that  the coordinates $\theta, \phi$ and the functions $s_k$
do not vary as the masses are changes.  
 The metric $ds^2_m$ when expressed in our coordinates is mass-dependent and
is given by 
$$ds^2 _m = \lambda^2 ds^2_1 \hskip 1cm (7.3a) $$
where
$$ds^2_1 ={1 \over 4}(d \phi ^2 + \cos(\phi)^2 d \theta^2) \hskip 1cm (7.3b)$$
(see [MontI], eq (5.6)  and Prop. 2) is the shape sphere metric when all the
masses are $1$,  and where the conformal factor
$\lambda$ is given  by 
$$\lambda = d(m)I_1/I ; \hskip .2cm d(m) = \sqrt{ 3 m_1 m_2 m_3/ M}  \hskip 1cm (7.3c)$$
(See eq (5.7) of [MontI]).
The metric we are to work with is thus 
$$ds^2 _J  = \tilde U\lambda^2 ds^2_1 \hskip 1cm (7.4).$$
Its curvature is given by lemma 4.1 :
$$\bar K = {{ 1} \over { \tilde U  \lambda^2}} \{ 4 - {1 \over 2} \Delta \log (\tilde U \lambda^2) \} \hskip 1cm (7.5).$$
where the
Laplacian $\Delta$ is with respect to the metric
$ds^2_1$.   The conformal factor is 
$$\eqalign{
\tilde U \lambda^2 & = d(m)^2 I U I_1^2/I^2 \cr
& = d(m) I_1 U (d(m) I_1/I) \cr
& = d(m) \hat U \lambda
}$$
where 
$$\hat U = I_1 U.$$
Then $$\Delta \log \tilde U \lambda^2 = \Delta \log \hat U + \Delta \log \lambda.$$
Set 
$$\hat s_i = s_i/I_1 = (1 - \gamma_k (\theta) \cos(\phi)  $$
so that
$$\hat U = \Sigma p_i/ \hat s_i $$
Since $\hat s_i$ is equal to the variable
$s_i$ as given by eq. (5.1.1b), and since the Laplacian
is the Laplacian of that section, we can continue
to   compute  as per section 5.
We have
$$ \hat U \lambda \bar K =  4 - {1 \over 2} \Delta \log (\hat U) - {1 \over 2} \Delta \log (\lambda) .$$

To compute the last term $- {1 \over 2} \Delta \log (\lambda)$ of (7.5)
we use the fact that  the metrics   $ds^2_m$ are 
$ds^2_1$  through $\lambda^2$, and that 
both have  curvature $4$. From lemma 4.1 it follows that
$$4 = {1 \over \lambda^2} (4 -  {1 \over 2} \Delta \log \lambda^2 \}$$
or
$$4\lambda^2 - 4 = - \Delta \log \lambda.$$
Equation (7.5) can then be rewritten 
$$ -\hat U^3 \lambda \bar K =  -\hat U^2(4 - {1 \over 2} \Delta \log (\hat U))
 + \hat U^2 (2 - 2 \lambda^2)  \hskip 1cm
(7.6).$$

{\bf Computation of $\Delta \hat U$.}
 To ease  notation, we drop the hats for this subsection,
so that   $U$ means the function $\hat U$ and
the   $s_k$ means the function 
$\hat s_k$. 
The   structural form of all the formulae and calculations of section 5
remains intact provided we  insert   the weightings
$p_i$ in the correct places.  

We proceed to the computation of $\Delta \log (\hat U))$.  Equations (4.2)
and (5.2.1) continue to hold with $\hat U$ in place of $U$, and equations
(5.2.2a,b) hold with $\hat s_i$ in place of $s_i$.  The analogues of
(5.2.3), (5.2.4) are
$$
{1 \over c}
\dphi(c \dphi U)   = ((s^2 - c^2)/c)\Sigma p_i \gamma_i /s_i ^2 
 + 2 s^2 \Sigma p_i\gamma_i ^2/ s_i^3 \hskip 1cm
\hskip 1cm (7.7a)$$
And
$$ {1 \over c^2} \dtheta \dtheta U 
  = {1 \over c} \Sigma p_i \gamma_i {\prime \prime}/s_i ^2
+ 2 \Sigma p_i (\gamma_i ^{\prime})^2/s_i ^3
\hskip 1cm (7.7b)$$
Now note that  algebraic steps going from (5.2.5) to (5.2.7) 
apply term by term, so can be carried over verbatim except that
the $i$th term must be multiplied by $p_i$. Thus (reverting to hats)
$$\Delta \hat U = 8 \hat U_4 \hskip 1cm (7.8a)$$
where $$\hat U_4 = \Sigma p_i/\hat s_i ^2 \hskip 1cm (7.8b).$$

The first line of  (5.3.2) becomes 
$$\|\nabla U \|^2  =  4(\Sigma s p_i\gamma_i/s_i^2)^2 +  4 (\Sigma p_i\gamma_i ^{\prime}/s_i ^2 )^2$$
The algebra which follows is essentially the same, leading to
$$\| \nabla U \|^2 = 4 S \hskip 1cm (7.9a)$$
where
$$S = 2 \Sigma p_i ^2/ s_i ^3 - \Sigma p_i ^2/s_i ^3 
-{3/2} \Sigma^{\prime} p_i p_j /s_i ^2 s_j ^2 +  2 \Sigma^{\prime} p_i p_j/s_i  s_j ^2
- \Sigma^{\prime} p_i p_j /s_i  s_j \hskip 1cm (7.9b).$$
Combining (7.8) and (7.9)  according to (4.2) (see also the steps (4.7)-(4.9)) yields the formula: 
$$- U^2 (4 - {1 \over 2} \Delta \log U)  = 3 \Sigma^{\prime} p_i p_j /s_i ^2 s_j ^2 - 2(\Sigma p_i /s_i)^2 
\hskip 1cm (7.10).$$

\vskip 1cm 
From earlier, we have   
$$ U^2 (2 - 2 \lambda^2) = 
 2 (\Sigma p_i /s_i)^2 - 2(\Sigma p_i /s_i)^2 d(m)^2 ({1 \over 3} \Sigma
s_i)^2/(\Sigma p_i s_i /M)^2, \hskip 1cm (7.11).$$
(We continue to use $s_i$ in place of $\hat s_i$.) Upon adding (7.10) and (7.11) there
is a   cancellation yielding:  
$$ - U^3 \lambda^2 \bar K =  3 \Sigma^{\prime} p_i p_j /s_i ^2 s_j ^2   - 2(\Sigma p_i /s_i)^2 d(m)^2 ({1 \over 3} \Sigma
s_i)^2/(\Sigma p_i s_i /M)^2.$$
A computation shows that
$$d(m)^2 M^2 = {3 \over 2}\Sigma^{\prime} p_i p_j$$
Recall that the $s_i$ (the previous $\hat s_i$'s) satisfy
${1 \over 3} \Sigma s_i =1$.  We finally obtain:
$$ - U^3 \lambda^2 \bar K =  3 \{ \Sigma^{\prime} p_i p_j /s_i ^2 s_j ^2  
 - \Sigma^{\prime} p_i p_j {{ (\Sigma p_i /s_i)^2 } \over {(\Sigma p_i s_i )^2}} \} 
\hskip 1cm (7.12).$$
Consequently 
$$\kappa = \sqrt{\Sigma^{\prime} {{p_i p_j} \over{ s_i ^2 s_j ^2}}}  
{{ (\Sigma p_i s_i)} \over {(\Sigma p_i /s_i )}} - \sqrt{  \Sigma^{\prime} p_i p_j }
\hskip 1cm (7.13)$$
governs the sign of the curvature,
with the curvature   $\bar K$   negative if $\kappa$
is positive, positive if $\kappa$ is negative, and zero
if $\kappa$ is zero.  We note that 
$\kappa$, and hence the curvature is zero 
at the Lagrange point $s_1 = s_2 =s_3 =1$,
and that this is true for all choices of the masses $p_i$.

\proclaim Theorem 6.  
For a Zariski-dense set of  mass distributions,
the sign of the curvature changes in a neighborhood of the
Lagrange point.

Proof. It suffices
to show that for a Zariski-dense set of  mass distributions
the differential  $d \kappa  \ne 0$ at the Lagrange point.
A differential form $\Sigma a_i ds_i$
represents zero on the shape sphere if and only if it is proportional to
$\Sigma ds_i$, the latter being the differential of the constraint
$\Sigma s_i = 3$ satisfied by the $s_i$ (which are the old $\hat s_i$'s). 
A computation shows
that  at  the Lagrange point $s_i =1$ we have  
$$\beta d \kappa = p_1 (p_2 ^2 + p_3 ^2) ds_1 + p_2 (p_1 ^2 + p_3 ^2) ds_2 + p_3 (p_1 ^2 + p_2 ^2) ds_3
\hbox{ mod }  \Sigma ds_i$$
where $\beta$ is a nonzero constant.
We thus want to know whether or not
the equality
$$(p_1 (p_2 ^2 + p_3 ^2), p_2 (p_1 ^2 + p_3 ^2) , p_3 (p_1 ^2 + p_2 ^2)) 
= (\lambda, \lambda, \lambda) \hskip 1cm (**) $$
can be satisfied for some $\lambda$.   
The right hand side of equation (**),
being homogeneous of degree $3$, 
  defines a polynomial map $\RP^2 \to \RP^2$
and we want to know if it is equal to the
constant map   $[1,1,1]$.  Because the map is  polynomial,
if we can exhibit a single point where the  inequality fails
then it must fail on a Zariski-dense set.    Plugging in
$p_1 = p_2 = 1, p_3 = a$ yields
$(p_1 (p_2 ^2 + p_3 ^2), p_2 (p_1 ^2 + p_3 ^2) , p_3 (p_1 ^2 + p_2 ^2))
= (1 + a^2, 1 + a^2, 2a)$
which is not proportional to $(1,1,1)$ unless $a = 1$.

QED 

\vskip .3cm 
{\bf 
Appendix A}

We prove
\proclaim Theorem A.  The set of initial conditions
within $H = 0,  I = 1 , c = 0$
whose solutions tend to a binary collision
of type $ij$ has nonempty interior. This fact
holds for all positive mass distributions. 

{\bf Proof of theorem A.}We use Newtonian time and Jacobi
coordinates For notational simplicity,  take $ij = 12$.
The   Jacobi coordinates are
$\zeta_1 = x_1 - x_2$, $\zeta_2 = x_3 - (m_1 x_1 +m_2 x_2)/(m_1 + m_2)$.
The distance to binary collision is
$$r = | \zeta_1| \hskip 1cm (1A).$$
We will exhibit a nonempty open set of
initial conditions at time $t =0$
 for which $r(t) = 0$
for some time $t < O(r(0))$.

The Hamiltonian is 
$$H = {1 \over 2} ( \mu_1 | \dot \zeta_1|^2 + 
\mu_2 | \dot \zeta_2 |^2)  -{{m_1 m_2} \over {r^2}} 
- W(\zeta_1, \zeta_2) \hskip 1cm (2A)$$
where $\mu_1 = m_1 m_2 /(m_1 + m_2)$,
$\mu_2 = m_3(m_1 + m_2)/M$, 
$W = m_1 m_3 /s_2 + m_2 m_3/s_1 $
and the squared distances $s_2$, $s_3$  can be expressed 
$|\zeta_2 + a_i
\zeta_1 |^2$ in Jacobi coordinates,  with
 mass-dependent
nonzero constants
$a_i$.   
The interaction term  $W$  satisfies the estimates  
$$|W| \le C_1 + C_2 \epsilon ,\hskip .3cm \hbox{ for } r < \epsilon \hskip 1cm 
(3A), $$ 
$$|{{\partial W} \over {\partial
\zeta_1}}|
\le C_1 + C_2 \epsilon 
\hskip .1cm ,  \hskip .1cm  |{{\partial W} \over {\partial
\zeta_2}}| \le C_1 + C_2 \epsilon , \hskip .3cm 
\hbox{ for  } r  \le \epsilon \hskip 1cm (4A)$$
($\epsilon$ sufficiently small), where $C_1, C_2$
are constants depending only on the masses.   

The equations of motion are 
$$\ddot \zeta_1 = - (m_1 + m_2){{\zeta_1} \over {r^4}}  + {1
\over
\mu_1}{{\partial W} \over {\partial \zeta_1}} \hskip 1cm (5A)$$
and
$$\ddot \zeta_2 =   {1 \over
\mu_2}{{\partial W} \over {\partial \zeta_2}} \hskip 1cm (6A)$$
Write 
$$J_1 = \zeta_1 \wedge \dot \zeta_1 \hskip 1cm (7A)$$
for the angular momentum (up to a factor of $\mu_1$)
of the 12 system.  We compute 
$$\eqalign{\dot J_1 & = \zeta_1 \wedge \ddot \zeta_1 \cr
&= \zeta_1 \wedge{1 \over
\mu_1}{{\partial W} \over {\partial \zeta_1}}}
$$
so that
$$|\dot J_1| \le Cr \hskip 1cm (8A).$$

Because $| \dot \zeta_1|^2 = \dot r^2 + J_1^2/ r^2$ we
have that 
$$r^2 H = {1 \over 2} (\mu_1 r^2 \dot r^2 + \mu_1 J_1^2 )
- m_1 m_2  + r^2 ( {1 \over 2} \mu_2 | \dot \zeta_2 |^2
- W ) \hskip 1cm (9A) $$

Now let $\zeta(t) = (\zeta_1 (t), \zeta_2 (t))$
be a solution satisfying the initial conditions
$r(0) < \epsilon,  H = 0, I =1, J = 0$.
From (6A) and  (4A)
we have that 
$$| \dot \zeta_2 (t)|^2 \le | \dot \zeta_2 (0)|^2 + Ct \hskip 1cm (10A),$$
for $t = O(1)$, provided $r(0) < \epsilon$.
Here $C$ depends only on the masses and $\epsilon$.  
 Letting $r \to 0$, we see from (9A) that if our solution
is to have a collision then
we must have
$$\lim r^2 \dot r^2 + \lim J_1 ^2 - 2 (m_1 + m_2) = 0,  \hskip 1cm (11A)$$
 where we have used $m_1 m_2 /\mu_1 = m_1 + m_2$.
But $r^2 \dot r^2 \ge 0$, so we must have
$$ 2 (m_1 + m_2) -\lim J_1 ^2  \ge 0 \hskip 1cm (12A)$$
 
We  argue   in the reverse. Suppose that   $ 2 (m_1 + m_2) - J_1 (0) ^2 $ is
sufficiently positive at the initial time $t=0$, and that $\dot r (0) < 0$ 
then (11A) forces $r^2 \dot r^2$ to be positive over
a finite  time interval. We will show that, upon  integration, 
this  will  force
$r(t) = 0$ in some finite time  $t = O(\sqrt{r(0)})$.
Note from the bounds (10A)   and the fact that $H= 0$
we have
$$|\mu_1  r^2 \dot r^2 + \mu_1 J_1 (0)^2 )
- 2 m_1 m_2  |  \le  K r(0) \hskip 1cm (13A) $$
for $0 \le t \le 1$ and for as long as $\dot r(t) < 0$. 
Here the constant $K$ depends only on the masses
and $\dot \zeta_2 (0)$.   Dividing by $\mu_1$ and using
$m_1 m_2 /\mu_1 = m_1 + m_2$ we arrive at 
$$|  r^2 \dot r^2 + J_1 (0)^2 )
- 2(m_1  + m_1 ) |  \le  K^* r(0)  \hskip 1cm (14A) $$ 
where $K^* = K/\mu_1$. 

We now impose the 
open condition
$$2(m_1 +  m_2)  - J_1 (0)^2 - K^* r(0) > \delta^2 \hskip 1cm 
(15A)
$$
on our initial conditions.  This will be the open condition
of theorem A.  
The positive constant    $\delta$ will
be constrained further
below.     
It follows from (15A) and (14A) that 
$$\delta^2  <  2(m_1 +  m_2) - J_1 (0)^2) - K^* r(0)  \le  r^2 \dot
r^2  \hskip 1cm 
(16A)$$
(16A) together with $\dot r(0) < 0$  forces   $\dot r < 0$ throughout the time interval in
question.   Thus $ - r \dot r > 0 $, and so we can take
square roots of  inequality  (16A) to obtain
$$\delta \le  - r \dot r  \hskip 1cm 
(17A) .$$  Taking negatives and integrating we find
that
$ -  \delta t \ge {1 \over 2} r(t)^2 - {1 \over 2}
r(0)^2$ or 
$$r(0)^2 - 2  \delta t  \ge r(t)^2 \hskip 1cm 
(18A).$$
This forces $r(t) = 0$ for some time
$t$ with $t \le r(0)^2 / 2 \delta $. 
In order that the collision time $t$ is $o(1)$ 
it is sufficient to take $\delta = O(r(0))$.   

We have proved that a 12 collision occurs
within a time $t = r(0)^2/ 2 \delta $ 
for all initial conditions satisfying (15A), 
$\dot r(0) < 0$,  and $r(0) < \epsilon$, 
where $\epsilon$ is small enough so that the inequalities  (3A, 4A) are in force.
This set of initial conditions is clearly open.
It remains to show that this set is nonempty.
Consider  the {\it collinear}  solution
having  $H =
0 = J$ and $I =1$.  (There are precisely two
such solutions, up to   time translation
and rotation, one for each arc of the equator which
ends in the 12 collision.)   These solutions  satisfy $J_1 = 0 $.
   In this case
(15A) reads
$m_1 + m_2  > K r(0) + \delta$.
and so will hold for $r (0)$ small provided
only that $\delta < m_1 + m_2$.
Since the solution tend to collision
it follows that (15A) is eventually in force
along the collinear solution, and hence that
our set of of initial conditions is nonempty.  

QED

\vskip 1cm

{\bf Acknowledgements.}  
I   dedicate this paper to the memory of my father.
I    acknowledge useful correspondences
with  Toshiaki Fujiwara, Alain Chenciner, Alain Albouy,
and   conversations with  Anatole Katok, Rafe Mazzeo, and with 
  Jeff Xia for pointing out that theorem 1
combined with an earlier
version of theorems 2, 3 and 4 ought to imply theorem 5.  

\medskip

{\bf References}

[AbMar] R. Abraham and J. Marsden, {\bf Foundations of Mechanics}, 
Benjamin-Cummings, [1978]. 

[AlbCh]  A. Albouy and A. Chenciner, {\it Le probl\'eme des $n$ corps et les
distances mutuelles}, Invent Math. 131 (1998), no. 1, 151--184.
 
[Arn]  V.I. Arnol'd, {\bf Mathematical Methods of Classical Mechanics},
Springer-Verlag, [1989]. 

[Ban] T. Banachiewitz,
{\it Sur un cas particulier du probleme
des trois corps}, CRAS, Paris, 142, (1906), pp 510-512.

[CGMS] A. Chenciner, J. Gerver, R. Montgomery R. and C. Sim\'o{\it Simple
choreographies of $N$ bodies: a preliminary study}
 in {\bf Geometry, Mechanics
and Dynamics}, 287--308, Springer,New York, 2002.

[ChMont] Chenciner A. and Montgomery R. {\it A remarkable periodic solution of the three-body
problem in the case of equal masses, Annals of Math., 152, pp. 881-901 (2000)}

[FerrTerr] Davide Ferrario and Susanna Terracini  
{\it  On the Existence of Collisionless Equivariant Minimizers for the Classical
n-body Problem.} Math ArXivs,         math-ph 0302022 [2003].

[Gordon] W. B. Gordon, {\it A minimizing property of keplerian orbits}, American
Journal of Mathematics, vol. 99,  $n^05$,   961-971, (1970).

[Fuji]  T. Fujiwara, H. Fukuda,  A. Kameyama,  H. Ozaki, M. Yamada,
{\it  Synchronised Similar Triangles for Three-Body Orbit
with Zero Angular Momentum}, arxiv.org/abs/math-ph/0404056.

[Hermann] R. Hermann, {\it On the differential geometry of foliations}, Ann. of
Math. (2), (1959), 445-457 

[Math] Mathematical Society of Japan, {\bf Encyclopedic Dictionary of
Mathematics},  by the Mathematical Society of Japan, ed. by S Iyanga and Y Kawada,
translated by K. O. May, The MIT Press, Cambridge, Massachussets, and
London, England, [1977]

[MontI]Richard Montgomery, {\it Infinitely Many Syzygies}, 
Archives for  Rational Mechanics and Analysis, 
v. 164 (2002), no. 4, 311--340, 2002.

[MontN] Richard Montgomery,  {\it The N-body problem, the braid group, and
action-minimizing periodic orbits}, Nonlinearity, vol. 11, no. 2, 363-376, 1998.

[MontR] Richard Montgomery, 
{\it  Geometric Phase of the Three-Body Problem}, Nonlinearity,
vol. 9, no. 5, 1341-1360, 1996.

[Moore] Cris Moore,  {\it Braids in Classical Gravity},  Physical Review
Letters 70, pp. 3675--3679, (1993).

[Morse] H. M. Morse, {\it A one-to-one representation of geodesics on a surface of
negative curvature}, Am. J. Math, {\bf 43}, no. 1, 33-51, 1921.

[Mumf] David Mumford, Caroline Series,David Wright, David  
{\bf Indra's pearls.} 
The vision of Felix Klein. 
Cambridge University Press, New York, 2002.

[Poin] Poincar\'e, ~H. [1896], {\it Sur les solutions p\'eriodiques et le
principe de moindre action}, .
 {\it C.R.A.S. Paris} {\bf 123}, 915--918.

\medskip
 
\smallskip


\bigskip

\vfill\eject
\end

\end